\def\UseSection{%%
      \numberwithin{equation}{section}
	\theoremstyle{plain}% default theorem style 
      \newtheorem{theorem}    {Theorem}[section]
      \DefineTheorems % Use this to define other environments to be 
}

\def\DefineTheorems{%%

	\theoremstyle{definition}% ``defn'' theorem style 

	\theoremstyle{definition}% ``remark'' theorem style 

}

\newcommand{\Qed}{\qed \vspace{.5cm}}

\def\eq#1\en{\begin{equation}#1\end{equation}}  
\def\eqsplit#1\ensplit{
	\begin{equation}\begin{split}#1\end{split}\end{equation}
	}
\def\eqalign#1\enalign{
	\begin{align}#1\end{align}
	}
\def\eqmul#1\enmul{
	\begin{multline}#1\end{multline}
	}
 
\makeatletter
\newcommand{\labelcounter}[2]{{%
	\stepcounter{#1}%	First, increase the ``countC'' by one.
	\protected@write\@auxout{}%
	{\string\newlabel{#2}{{\csname the#1\endcsname}{\thepage}}}%
	{\ref{#2}}%	Finally, make sure to refer to this label, 
	}}
\makeatother
\newcommand{\Nbold} {{\mathbb N}}

\newcommand{\Zbold} {{\mathbb Z}}
\documentclass[12pt]{article}
\usepackage[psamsfonts]{amsfonts}
\usepackage{amsmath,amssymb,amsthm}

\usepackage{tikz,graphicx}
\usetikzlibrary{arrows}

\usepackage{hyperref}

\usepackage{epic,eepic}

\newtheorem{THM}{Theorem}[section]
\newtheorem{COR}[THM]{Corollary}

\newtheorem{REM}[THM]{Remark}
\newtheorem{LEM}[THM]{Lemma}

\newtheorem{DEF}[THM]{Definition}

\newcommand{\nn}{\nonumber}

\UseSection  %Necessary to define Numbering scheme for Theorem, etc.
\newcommand{\hlf}{\frac{1}{2}}
\newcommand{\ra}{\rightarrow}

\renewcommand{\to}      {\rightarrow}

\newcounter{countC}  % Defined the counter ``countC''
\newcounter{countR}  % Defined the counter ``countR''
\newcommand{\Z}{\Zbold}
\newcommand{\N}{\Nbold}

\newcommand{\cb}{\overset{b}{\rightarrow}}

\newcommand{\mc}[1]{\mathcal{#1}}
\newcommand{\mP}{\mathbb{P}}

\newcommand{\LRRL}{\begin{picture}(,)
\put(0,4){$\leftrightarrows$}
\put(0,-4){$\rightleftarrows$}
\end{picture}\hspace{.5cm}}

\newcommand{\RLLR}{\begin{picture}(,)
\put(0,4){$\rightleftarrows$}
\put(0,-4){$\leftrightarrows$}
\end{picture}\hspace{.5cm}}

\newcommand{\blank}[1]{}

\newcommand{\CL}{\preccurlyeq }

\newcommand{\TL}{\trianglelefteq }

\newcommand{\ACL}{\overset{\scriptscriptstyle\mc{A}}{\CL}}
\renewcommand{\ACL}{\CL^{\scriptscriptstyle\mc{A}}}

\title  {
      A combinatorial result with applications\\ to self-interacting random walks.
      }

\author{
Mark Holmes\footnote{Department of Statistics, University of Auckland.  E-mail {\tt
holmes@stat.auckland.ac.nz}} \and Thomas S. Salisbury \footnote{Department of Mathematics and Statistics, York University. E-mail {\tt
salt@yorku.ca}}}

\begin{document}

\maketitle

\begin{abstract}
We give a series of combinatorial results that can be obtained from any two collections (both indexed by $\Z\times \N$) of left and right pointing arrows that satisfy some natural relationship.   When applied to certain self-interacting random walk couplings, these allow us to reprove some known transience and recurrence results for some simple models. We also obtain new results for one-dimensional multi-excited random walks and for random walks in random environments in all dimensions.
\end{abstract}

%\tableofcontents
\section{Introduction}
\label{sec:intro}
Coupling is a powerful tool for proving certain kinds of properties of random variables or processes.  A coupling of two random processes $X$ and $Y$ typically refers to defining random variables $X'$ and $Y'$ {\em on a common probability space} such that $X'\sim X$ (i.e.~$X$ and $X'$ are identically distributed) and $Y'\sim Y$.  There can be many ways of doing this, but generally one wants to define the probability space such that the {\em joint distribution} of $(X',Y')$ has some property.  For example, suppose that $X=\{X_n\}_{n\ge 0}$ and $Y=\{Y_n\}_{n\ge 0}$ are two nearest-neighbour simple random walks in 1 dimension with drifts $\mu_X\le \mu_Y$ respectively.  One can define $X'\sim X$ and $Y'\sim Y$ on a common probability space so that $X'$ and $Y'$ are independent, but one can also define $X''\sim X$ and $Y''\sim Y$ on a common probability space so that $X''_n\le Y''_n$ for all $n$ with probability 1. 

Consider now a nearest-neighbour random walk $\{X_n\}_{n\ge 0}$ on $\Z^d$ that has transition probabilities $(2d)^{-1}$ of stepping in each of the $2d$ possible directions, except on the {\em first departure from each site}. On the first departure, these are also the transition probabilities for stepping to the left and right in any coordinate direction other than the first. But in the first coordinate, the transition probabilities are instead $(2d)^{-1}(1+\beta)$ (right) and $(2d)^{-1}(1-\beta)$ (left), for some fixed parameter $\beta\in [0,1]$.  This is known as an excited random walk \cite{BW03} and the behaviour of these and more general walks of this kind has been studied in some detail since 2003.  For this particular model, it is known \cite{BR07} that for $d\ge 2$ and $\beta>0$, there exists $v_{\beta}=(v^{[1]}_{\beta},0,\dots,0)\in \Z^d$ with $v^{[1]}_{\beta}>0$ such that $\lim_{n\ra \infty}n^{-1}X_n=v_{\beta}$ with probability 1.  When $d=1$ the model is recurrent (0 is visited infinitely often) except in the trivial case $\beta=1$.  It is plausible that  $v^{[1]}_{\beta}$ should be a non-decreasing function of $\beta$ (i.e.~increasing the local drift should increase the global drift) but this is not known in general.

A natural first attempt at trying to prove such a monotonicity result would be as follows:  given $0<\beta_1<\beta_2\le 1$, construct a coupling of excited random walks $X$ and $Y$ with parameters $\beta_1$ and $\beta_2>\beta_1$ respectively such that with probability 1, $X_n^{[1]}\le Y_n^{[1]}$ for all $n$.  Thus far no one has been able to construct such a coupling, and the monotonicity of $v^{[1]}_{\beta}$ as a function of $\beta$ remains an open problem in dimensions $2\le d\le 8$.  In dimensions $d\ge 9$ this result has been proved \cite{HH09mono} using a somewhat technical expansion method, as well as rigorous numerical bounds on simple random walk quantities.  More general models in 1 dimension have been studied, and some monotonicity results \cite{Zern05} have been obtained via probabilistic arguments but without coupling.  This raises the question of whether or not one can obtain proofs of these kinds of results using a coupling argument that has weaker aims e.g.~such that $\max_{m\le n}X_m^{[1]}\le \max_{m\le n}Y_m^{[1]}$ for all $n$, rather than $X_n^{[1]}\le Y_n^{[1]}$ for all $n$.

This paper addresses this issue in 1-dimension.  We study relationships between completely deterministic (non-random) 1-dimensional systems of arrows that may prove to be of independent interest in combinatorics.  Each system $\mc{L}$ of arrows defines a sequence $L$ of integers.  We show that under certain natural local conditions on arrow systems $\mc{L}$ and $\mc{R}$, one obtains relations between the corresponding sequences such as $\max_{m\le n}L_m^{[1]}\le \max_{m\le n}R_m^{[1]}$ for all $n$ (while it's still possible that $L_n^{[1]}> R_n^{[1]}$ for some $n$).   

These may be applied to certain random systems of arrows, to give self-interacting random walk couplings. Doing so, one can obtain results about the (now random) sequence $R_n$ {\em if} $L_n$ (also random) is well understood, and vice versa.  This yields alternative proofs of some existing results, as well as new non-trivial results about so-called multi-excited random walks in 1 dimension and some models of random walks in random environments in all dimensions -- see e.g.~\cite{HS_RWDRE}. To be a bit more precise, in \cite{HS_RWDRE} a projection argument applied to some models of random walks in random environments (in all dimensions) gives rise to a one-dimensional random walk $Y$, which can be coupled with a one-dimensional multi-excited random walk $Z$ (both walks depending on a parameter $p$) so that for every $j\in \Z$ and every $r\ge 1$:
\begin{itemize}
\item[(i)] If $Y$ goes left on its $r$th visit to $j$ then so does $Z$ (if such a visit occurs), and therefore
\item[(ii)] If $Z$ goes right on its $r$th visit to $j$ then so does $Y$ (if such a visit occurs).
\end{itemize}
Explicit conditions ($p>\frac{3}{4}$ in this case) governing when $Z_n\ra \infty$ as $n\ra \infty$ are given in \cite{Zern05}.  One would like to conclude that also $Y_n\ra \infty$ (whence the original random walk in $d$-dimensions returns to its starting point only finitely many times) when $p>\frac{3}{4}$.  This can be achieved by applying the result of this paper to the coupling mentioned above.  

The main contributions of this paper are: 
 combinatorial results concerning sequences defined by arrow systems satisfying certain natural local relationships (see Theorem \ref{thm:main}); some non-trivial counterintuitive examples; and application of these combinatorial results with non-monotone couplings to obtain new results in the theory of random walks.

\subsection{Arrow systems}
\label{sec:arrows}
A collection $\mc{E}=(\mc{E}(x,r))_{x \in \Z, r\in \mathbb{N}}$, where $\mc{E}(x,r)\in \{\leftarrow,\rightarrow\}$ is the arrow above the vertex $x\in \Z$ at level $r\in \N$,  is called an {\em arrow system}.  This should be thought of as an infinite (ordered) stack of arrows rising above each vertex in $\Z$.

In a given arrow system $\mc{E}$, let $\mc{E}_{\leftarrow}(j,r)$ denote the number of $\leftarrow$ arrows, out of the first $r$ arrows above $j$.  As $r$ increases, this quantity counts the number of $\leftarrow$'s appearing in the arrow columns above $j$.  Similarly define $\mc{E}_{\rightarrow}(j,r)=r-\mc{E}_{\leftarrow}(j,r)$.  We can define a sequence $E=\{E_n\}_{n\ge 0}$ by setting $E_0=0$ and letting $E$ evolve by taking one step to the left or right (at unit times),  according to the lowest arrow of the $\mc{E}$-stack at its current location, and then deleting that arrow.  In other words, if $\#\{0\le m\le n:E_m=E_n\}=k$ then $E_{n+1}=E_n+1$ if $\mc{E}(E_n,k)=\rightarrow$ (resp.~$E_{n+1}=E_n-1$ if $\mc{E}(E_n,k)=\leftarrow$).  

%In Section \ref{sec:excursion} we'll briefly discuss the connection between these sequences and arrow systems and the theory of excursions and trees.

\begin{DEF}[$\mc{L}\CL\mc{R}$]
\label{def:cl}
Given two arrow systems $\mc{L}$ and $\mc{R}$, we write $\mc{L}\CL\mc{R}$ if for each $j\in \Z$ and each $r\in \mathbb{N}$, 
\[\mc{L}_{\leftarrow}(j,r)\ge \mc{R}_{\leftarrow}(j,r) \qquad (\text{and hence also } \mc{L}_{\rightarrow}(j,r)\le \mc{R}_{\rightarrow}(j,r)).\]
\end{DEF}

\begin{DEF}[$\mc{L}\TL\mc{R}$]
\label{def:tl}
We write $\mc{L}\TL\mc{R}$ if for each $j\in \Z$ and each $r\in \mathbb{N}$,
\[\mc{L}(j,r)=\rightarrow \quad \Rightarrow \mc{R}(j,r)=\rightarrow.\]
%there is no $x\in \Z$, $r\in \mathbb{N}$ such that both $\mc{L}_{\rightarrow}(j,r)=\mc{L}_{\rightarrow}(j,r-1)+1$ and $\mc{R}_{\rightarrow}(j,r)=\mc{R}_{\rightarrow}(j,r-1)$.  This simply says that there is no vertex and level at which the $\mc{L}$-arrow at that level above that vertex points right, while the corresponding arrow in $\mc{R}$ points left.  
\end{DEF}
\noindent It is easy to see that $\mc{L}\TL \mc{R}$ implies $\mc{L}\CL \mc{R}$.

Now define two paths/sequences $\{L_n\}_{n\ge 0}$ and $\{R_n\}_{n\ge 0}$ in $\Z$ according to the arrows in $\mc{L}$ and $\mc{R}$ respectively as above (in particular $L_0=R_0=0$).  Since each arrow system determines a unique sequence, but a given sequence may be obtained from multiple different arrow systems, we write $L\CL R$ (resp.~$L\TL R$) if there exist $\mc{L}\CL \mc{R}$ (resp.~$\mc{L}\TL \mc{R}$) whose corresponding sequences are $L$ and $R$ respectively.  Note that when $\mc{L}\TL\mc{R}$, the paths $Z=L$ and $Y=R$ constructed from $\mc{L}$ and $\mc{R}$ as above automatically satisfy the conditions (i) and (ii) appearing at the beginning of Section \ref{sec:intro}.  

An arrow system $\mc{E}$ is said to be {\em 0-right recurrent} if in the new system $\mc{E}_+$ defined by $\mc{E}_+(0,i)=\rightarrow$ for all $i\ge 1$, and $\mc{E}_+(x,i)=\mc{E}(x,i)$ for all $i\ge 1$ and $x>0$, $E_{+,n}=0$ infinitely often.  

%%%%%%%%%%%%%%%%
%Similarly $\mc{E}$ is {\em o-left recurrent} if in the system $\mc{E}_-$ defined by $\mc{E}_-(0,i)=\leftarrow$ for %all $i\ge 1$
%% and  $x\le 0$
%, and $\mc{E}_-(x,i)=\mc{E}(x,i)$ for all $i\ge 1$ and $x<0$, $E_{-,n}=o$ infinitely often.

%\begin{DEF}[$|\mc{I}|\CL\mc{E}$]
%\label{def:abs_cl}
%We write $|\mc{I}|\CL\mc{E}$ if:
%\begin{itemize}
%\item  for each $j>0$ and $r\ge1$, $\mc{I}_{\leftarrow}(j,r)\ge \mc{E}_{\leftarrow}(j,r)$, and 
%\item  for each $j<0$ and $r\ge1$, $\mc{I}_{\rightarrow}(j,r) \ge \mc{E}_{\rightarrow}(j,r)$.
%\end{itemize}
%That is, if $\mc{I}$ has higher counts for arrows pointing toward $o$.
%\end{DEF}

The main result of this paper is the following theorem, in which $n_{E,t}(x)=\#\{k\le t:E_k=x\}$ (see also Corollary \ref{cor:transience} in the case that $L$ is transient to the right).
\begin{THM}
\label{thm:main}
Suppose that $\mc{L}\CL \mc{R}$.  Then 
\begin{itemize}
\item[(i)] $\liminf_{n\ra \infty} L_n\le \liminf_{n\ra \infty} R_n$;
\item[(ii)] $\limsup_{n\ra \infty} L_n\le \limsup_{n\ra \infty} R_n$;
\item[(iii)] Let $a_n\le n$ be any increasing sequence, with $a_n\to\infty$. If there exists $x \in \Z$ such that $R\ge x$ infinitely often then $\limsup_{n\ra \infty} \frac{L_n}{a_n}\le \limsup_{n\ra \infty} \frac{R_n}{a_n}$.
\item[(iv)] If $n_{R,t}(x)>n_{L,t}(x)$ then $n_{R,t}(y)\ge n_{L,t}(y)$ for every $y>x$.
\item[(v)] If $\mc{R}$ is $0$-right recurrent then so is $\mc{L}$.
\end{itemize}
\end{THM}

%\begin{THM}
%\label{thm:recurrent}
%Suppose that $\mc{I}\CL\mc{E}$.  Then if $\mc{E}$ is $o$-right recurrent then so is $\mc{I}$.
%%\item[(ii)] If $\mc{E}$ is $o$-left recurrent then so is $\mc{I}$.  
%\end{THM}

As $\frac{L_n}{n}$ represents the average speed of the sequence $L$, up to time $n$, in many applications the sequence of interest in Theorem \ref{thm:main} (iii) will be $a_n=n$.  Part (ii) of Theorem \ref{thm:main} actually follows from part (i) by a simple mirror symmetry argument.  There is a symmetric version of (iii), but one must be careful.  Part (iii) obviously implies that if $u=\lim n^{-1}R_n$ and $l=\lim n^{-1}L_n$ both exist then $l\le u$, however we show in Section \ref{sec:counterex1} that $L\TL R$ does not imply that $\liminf\frac{L_n}{n}\le \liminf\frac{R_n}{n}$.  The mirror image (about 0) of the counterexample in Section \ref{sec:counterex1} also shows that (iii) is not true in general if we drop the condition that $L\ge x$ infinitely often, for some $x$.  
One might also conjecture that if $L\TL R$ then the amount of time that $R>L$ is at least as large as the amount of time that $R<L$.  This is also false as per a counterexample in Section \ref{sec:counterex2}.

The remainder of the paper is organised as follows.  Section \ref{sec:basic} contains the basic combinatorial relations which are satisfied by the arrow systems and their corresponding sequences.  These will be needed in order to prove our first results. Section \ref{sec:results} gives various consequences of the relationship $\mc{L}\CL \mc{R}$ between two arrow systems, and includes the proofs of the main results of the paper.  Section \ref{sec:counterex} contains the counterexamples described above.  
%Section \ref{sec:excursion} briefly discusses the relationship between the existing theory of excursions and trees, and our arrow systems.  
Finally Section \ref{sec:applications} contains applications of our results in the study of self-interacting random walks.

\section{Basic relations}
\label{sec:basic}
Given an arrow system $\mc{E}$ and $t\ge 0$, let $n_{E,t}(x)=\#\{k\le t:E_k=x\}$ and $n_{E,t}(x,y)=\#\{k\le t:E_{k-1}=x,E_k=y\}$.  Then the following relationships hold:

\eqalign
n_{E,t}(x)&=\delta_{x,0}+n_{E,t}(x-1,x)+n_{E,t}(x+1,x)\label{eq:fact1}\\
n_{E,t}(x)&=\delta_{E_t,x}+n_{E,t}(x,x+1)+n_{E,t}(x,x-1)\label{eq:fact2}\\
t+1&=\sum_{i=-\infty}^{\infty}n_{E,t}(i).\label{eq:fact3}
\enalign
Relation (\ref{eq:fact1}) says that every visit to $x$ is either from the left or right, except for the first visit if $x=0$.  Relation (\ref{eq:fact2}) is similar, but in terms of departures from $x$.  The sum in (\ref{eq:fact3}) is in fact a finite sum since $n_{E,t}(i)=0$ for $|i|>t$. 

Next 
\eqalign
n_{E,t}(x,x+1)&=\mc{E}_{\rightarrow}(x,n_{E,t}(x)-I_{E_t=x})\label{eq:R1}\\
n_{E,t}(x,x-1)&=\mc{E}_{\leftarrow}(x,n_{E,t}(x)-I_{E_t=x})\label{eq:R2},
\enalign
where e.g.~relation (\ref{eq:R1}) says that the number of departures from $x$ to the right is the number of ``used'' right arrows at $x$.

Finally, 
\eqalign
n_{E,t}(x,x+1)+I_{x+1\le 0}I_{E_t\le x}&=n_{E,t}(x+1,x)+I_{x\ge 0}I_{E_t\ge x+1},\label{eq:pm1}
\enalign
which says that the number of moves from $x$ to $x+1$ is closely related to the number of moves from $x+1$ to $x$.  They may differ by 1 depending on the position of $x$ relative to 0 and the current value of the sequence.  For example, if $0\le x<E_t$ then the number of moves from $x$ to $x+1$ up to time $t$ is one more than the number of moves from $x+1$ to $x$ up to time $t$.
%If $L\CL R$ then for all $l,r$,
%\eqalign
%n_{R,r}(x-1)\ge n_{L,l}(x-1) &\Rightarrow n_{R,r}(x-1,x)\ge n_{L,l}(x-1,x) \label{eq:imp1}\\
%n_{R,r}(x+1)\le n_{L,l}(x+1) &\Rightarrow n_{R,r}(x+1,x)\le n_{L,l}(x+1,x) \label{eq:imp2}.
%\enalign
%To see (\ref{e:imp1}), for example, observe that 
%\begin{align*}
%n_{R,r}(x-1,x)&=\mc{R}_{\rightarrow}(x-1,n_{R,r}(x-1))
%\ge\mc{R}_{\rightarrow}(x-1,n_{L,l}(x-1))\\
%&\ge\mc{L}_{\rightarrow}(x-1,n_{L,l}(x-1))
%=n_{L,l}(x-1,x).
%\end{align*}

\section{Implications of $\mc{L}\CL \mc{R}$.}
\label{sec:results}
In this section we always assume that $\mc{L}\CL\mc{R}$.  The results typically have symmetric versions using the fact that $\mc{L}\CL \mc{R} \iff -\mc{R}\CL -\mc{L}$, which is equivalent to considering arrow systems reflected about 0.  We divide the section into two subsections based roughly on the nature of the results and their proofs.  

For $x \in \Z$ and $k\ge 0$, let $T_{L}(x,k)=\inf\{t\ge 0:n_{L,t}(x)=k\}$, and $T_{R}(x,k)=\inf\{t\ge 0:n_{R,t}(x)=k\}$.

\subsection{Results obtained from the basic relations}
\label{sec:appl_basic}
The proofs in this section are based on applications of the basic relations of Section \ref{sec:basic}.  The first few results are somewhat technical, but will be used in turn to prove some of the more appealing results.  Roughly speaking they describe how the relative numbers of visits of $L$ and $R$ to neighbouring sites $x-1$ and $x$ relate to each other.  
\begin{LEM}
\label{lem:exists}
If $L$ hits $x$ at least $k\ge 1$ times and $R$ is eventually to the left of $x$ after fewer than $k$ visits to $x$, then there exists a site $y<x$ that $R$ hits at least $n_{L,T_L(x,k)}(y)$ times.
\end{LEM}
\proof Fix $x,k$ and let $T=T_L(x,k)$ and $y_0:=\inf\{z\le x:n_{L,T}(z)>0\}\le 0$.  If $y_0=x$ then the first $k-1$ arrows at $x$ are all right arrows, i.e.~$\mc{L}_{\rightarrow}(y_0,k-1)=k-1$.  Then also $\mc{R}_{\rightarrow}(y_0,k-1)=k-1$ so $R$ cannot be to the left of $x$ after fewer than $k$ visits.  Similarly if $y_0<x$ then the first $n_{L,T}(y_0)$ arrows at $y_0$ are all right arrows, i.e.~$\mc{L}_{\rightarrow}(y_0,n_{L,T}(y_0))=n_{L,T}(y_0)$, and so also $\mc{R}_{\rightarrow}(y_0,n_{L,T}(y_0))= n_{L,T}(y_0)$.  
%It follows that $y_0<x$. [This is obvious if $x_0>0$. While if $y_0=x_0\le 0$ we would have a contradiction to the hypothesis that $R$ is to the left of $x$ within $k$ visits.] 
Therefore either $R$ visits $y_0$ at least $n_{L,T}(y_0)$ times or it stays in $(y_0,x)$ infinitely often, whence it must visit some site $y \in (y_0,x)$ at least $n_{L,T}(y)$ times as required.\Qed

\vspace{.5cm}
Let $n_L(x)=n_{L,\infty}(x)$ and $n_R(x)=n_{R,\infty}(x)$.

\begin{LEM}
\label{lem:x-1,x}
If $R$ hits $x-1$ at least $n_L(x-1)$ times then either 
\begin{enumerate}
\item $n_R(x)\ge n_L(x)$, or
%$R$ hits $x$ at least $n_L(x)$ times, or 
\item $R$ is always to the right of $x$ after fewer than $n_L(x)$ visits. 
($\Rightarrow \liminf_{n\ra \infty}R_n>x$)
\end{enumerate}
\end{LEM}
\proof Assume that the first claim fails, so in particular $n_{R}(x)<\infty$.  Let $T=\inf\{t:n_{L,t}(x)=n_R(x)+1\}$.
Then $T<\infty$ so $L_T=x$. Choose $r$ sufficiently large so that $R_t\ne x$ for any $t\ge r$, $R_r\ne x-1$, and $n_{R,r}(x-1)\ge n_{L,T}(x-1)$.  
%Let $T=T_{L}(x,n_{R,r}(x)+1)=\inf\{l\ge 0:n_{L,l}(x)=n_{R,r}(x)+1\}$.  
Then by (\ref{eq:fact1}) applied to $L$ at time $T$, and also to $R$ at time $r$,
\eqalign
n_{R,r}(x)+1&= n_R(x)+1=n_{L,T}(x)=n_{L,T}(x-1,x)+n_{L,T}(x+1,x)+\delta_{0,x}\nn\\
n_{R,r}(x)&=\delta_{x,0}+n_{R,r}(x-1,x)+n_{R,r}(x+1,x).\nn
\enalign
Subtracting one from the other and rearranging we obtain
\eq
n_{R,r}(x-1,x)-n_{L,T}(x-1,x)+n_{R,r}(x+1,x)+1=n_{L,T}(x+1,x).\nn
\en
Now $n_{L,T}(x+1,x)=n_{L,T}(x,x+1)+I_{x+1\le 0}$ from (\ref{eq:pm1}), so
%and $n_{R}(x,x+1)+I_{x+1\le 0}I_{R_r\le x}&=n_{R}(x+1,x)+I_{x\ge 0}I_{R_r\ge x+1}$  Hence,
%\eq
%n_{R,r}(x-1,x)+n_{R,r}(x+1,x)+1-n_{L,T}(x-1,x)-[n_{L,T}(x,x+1)+I_{x+1\le 0}]=0.
%\en
%Rearranging we obtain
\eq
n_{R,r}(x+1,x)+1+[n_{R,r}(x-1,x)-n_{L,T}(x-1,x)]= n_{L,T}(x,x+1)+I_{x+1\le 0}.\label{eq:elephant}
\en
Using (\ref{eq:R1}) and the fact that $R_r\ne x$, then $\mc{L}\CL \mc{R}$, then the fact that $n_{L,T}(x)=1+n_{R,r}(x)$, and finally again using (\ref{eq:R1}) and the fact that $L_T=x$ we obtain
\eq
n_{R,r}(x,x+1)=\mc{R}_{\rightarrow}(x,n_{R,r}(x))\ge\mc{L}_{\rightarrow}(x,n_{R,r}(x))
=\mc{L}_{\rightarrow}(x,n_{L,T}(x)-1)=n_{L,T}(x,x+1).\nn
\en
Using this bound in (\ref{eq:elephant}) yields
\eq
n_{R,r}(x+1,x)+1+[n_{R,r}(x-1,x)-n_{L,T}(x-1,x)]\le n_{R,r}(x,x+1)+I_{x+1\le 0}.\label{eq:lineabove}
\en
Using the fact that $R_r\ne x-1$ and applying (\ref{eq:R1}) to $R_r$ at $x-1$, then using $n_{R,r}(x-1)\ge n_{L,T}(x-1)$, then $\mc{L}\CL \mc{R}$, and finally using the fact that $L_T\ne x-1$ and applying (\ref{eq:R1}) to $L_T$ at $x-1$, we have that 
\begin{align*}
n_{R,r}(x-1,x)
&=\mc{R}_{\rightarrow}(x-1,n_{R,r}(x-1))
\ge\mc{R}_{\rightarrow}(x-1,n_{L,T}(x-1))\\
&\ge\mc{L}_{\rightarrow}(x-1,n_{L,T}(x-1))
= n_{L,T}(x-1,x).
\end{align*}
Therefore by (\ref{eq:lineabove}),  and then (\ref{eq:pm1})
\eq
n_{R,r}(x+1,x)+1\le n_{R,r}(x,x+1)+I_{x+1\le 0}\le n_{R,r}(x+1,x)+I_{R_r\ge x+1}.\nn
\en
Therefore $R_r\ge x+1$, so in fact
$R_t>x$ for every $t\ge r$. Moreover $n_{R,r}(x)=n_R(x)<n_L(x)$, which shows (b). 
\Qed

\begin{LEM}
 \label{lem:main1}
Let $x\in \Z$, and suppose that for some $k>0$, $n_L(x)\ge k$ and $n_R(x)\ge k$.  Then $n_{R,T_R(x,k)}(x-1)\le n_{L,T_L(x,k)}(x-1)$.
\end{LEM}
\proof  Let $T=T_L(x,k)<\infty$ and $S=T_R(x,k)<\infty$.  Then $R_S=x>x-1$, so from (\ref{eq:pm1}) and (\ref{eq:R2})
\[n_{R,S}(x-1,x)=n_{R,S}(x,x-1)+I_{x\ge 1}=\mc{R}_{\leftarrow}(x,k-1)+I_{x\ge 1}.\]
Similarly 
\[n_{L,T}(x-1,x)=n_{L,T}(x,x-1)+I_{x\ge 1}=\mc{L}_{\leftarrow}(x,k-1)+I_{x\ge 1}.\]
Since $\mc{R}_{\leftarrow}(x,k-1)\le \mc{L}_{\leftarrow}(x,k-1)$ it follows that $n_{R,S}(x-1,x)\le n_{L,T}(x-1,x)$.  
Finally,
\[\mc{R}_{\rightarrow}(x-1,n_{R,S}(x-1))=n_{R,S}(x-1,x) \text{ and }n_{L,T}(x-1,x)=\mc{L}_{\rightarrow}(x-1,n_{L,T}(x-1))\]
whence $\mc{R}_{\rightarrow}(x-1,n_{R,S}(x-1))\le \mc{L}_{\rightarrow}(x-1,n_{L,T}(x-1))$.  Since the $n_{R,S}(x-1)$-th arrow at $x-1$ is $\rightarrow$ by definition of $S$ (and similarly for $n_{L,T}(x-1)$ and $T$) this implies that $n_{R,S}(x-1)\le n_{L,T}(x-1)$ as required.
\Qed

\begin{LEM}
If $T=T_L(x,k)<\infty$  and $R$ stays to the right of $x$ after fewer than $k$ visits to $x$ then $n_R(x-1)\le n_{L,T}(x-1)$.
\label{lem:main1.5}
\end{LEM}
\proof Assume that $n_R(x-1)>0$, otherwise there is nothing to prove. Let $S'=\sup\{t:R_t=x\}$.  Then $R_{S'}=x$, $\mc{R}(x-1,n_{R,S'}(x-1))=\rightarrow$ and $\mc{R}(x,n_{R,S'}(x))=\rightarrow$.  By (\ref{eq:pm1}) applied at $x-1$, and then using (\ref{eq:R2}), and  finally the fact that $\mc{R}(x,n_{R,S'}(x))=\rightarrow$, 
\[n_{R,S'}(x-1,x)=n_{R,S'}(x,x-1)+I_{x\ge 1}=\mc{R}_{\leftarrow}(x,n_{R,S'}(x)-1)+I_{x\ge 1}=\mc{R}_{\leftarrow}(x,n_{R,S'}(x))+I_{x\ge 1}.\]
Therefore by (\ref{eq:R1}),
\eqalign
&\mc{R}_{\rightarrow}(x-1,n_{R,S'}(x-1))=n_{R,S'}(x-1,x)=\mc{R}_{\leftarrow}(x,n_{R,S'}(x))+I_{x\ge 1}.\label{eq:platypus}
\enalign
Since $n_{R,S'}(x)<k=n_{L,T}(x)$ we have $\mc{R}_{\leftarrow}(x,n_{R,S'}(x))\le \mc{L}_{\leftarrow}(x,n_{L,T}(x)-1)$, therefore the right hand side of (\ref{eq:platypus}) is bounded above by
\eqalign
\mc{L}_{\leftarrow}(x,n_{L,T}(x)-1)+I_{x\ge 1}&=n_{L,T}(x,x-1)+I_{x\ge 1}\nn\\
&=n_{L,T}(x-1,x)=\mc{L}_{\rightarrow}(x-1,n_{L,T}(x-1)),\nn
\enalign
where we have used (\ref{eq:R2}), followed by (\ref{eq:pm1}), and then (\ref{eq:R1}).  We have shown that
\[\mc{R}_{\rightarrow}(x-1,n_{R,S'}(x-1))\le \mc{L}_{\rightarrow}(x-1,n_{L,T}(x-1)).\]
Since $\mc{R}(x-1,n_{R,S'}(x-1))=\rightarrow$, this implies that $n_{R,S'}(x-1)\le n_{L,T}(x-1)$ as required.
\Qed

\subsection{Results obtained by contradiction}
\label{sec:contra}
The results in this section include less technical results than those of the previous section.  Roughly speaking their proofs will be based on contradiction arguments that proceed as follows.  Suppose that we have already proved a statement $A$ whenever $\mc{L}\CL \mc{R}$.  We now want to prove a statement $B$ whenever $\mc{L}\CL \mc{R}$.  Assume that for some $\mc{L}$, $\mc{R}$ with $\mc{L}\CL \mc{R}$, $B$ is false.  Construct two new systems $\mc{L}'\CL \mc{R}'$ from $\mc{L}$ and $\mc{R}$ such that statement $A$ is violated for $\mc{L}'$ and $\mc{R}'$.  This gives a contradiction, hence there was no such example where $\mc{L}\CL \mc{R}$ but $B$ is false.
\begin{LEM}
\label{lem:main2}
Let $x\in \Z$, and suppose that $n_R(x)<k\le n_L(x)$.  Then $n_R(x-1)\le n_{L,T_L(x,k)}(x-1)$ and $\liminf R_n>x$ (i.e.~
%$L$ hits $x$ at least $k$ times.  
$R$ is forever to the right of $x$ after fewer than $k$ visits to $x$ and at most $n_{L,T_L(x,k)}(x-1)$ visits to $x-1$).
\end{LEM}
\proof By Lemma \ref{lem:main1.5}, it is sufficient to prove that under the hypotheses of the lemma, $R$ is to the right of $x$ infinitely often.    Suppose instead that $R$ is forever to the left of $x$ (after fewer than $k$ visits to $x$).  Then we may define two new systems $\mc{R}'$ and $\mc{L}'$ by forcing every arrow at $x$ at level $k$ and above to be $\rightarrow$.  To be precise, given an arrow system $\mc{E}$ we'll define $\mc{E}'$ by $\mc{E}'(y,\cdot)=\mc{E}(y,\cdot)$ for all $y\ne x$, $\mc{E}'(x,j)=\mc{E}(y,j)$ for all $j<k$, and $\mc{E}'(x,j)=\rightarrow$ for every $j\ge k$.  Clearly $\mc{L}'\CL\mc{R}'$ and $T'=T_{L'}(x,k)=T$.  The sequences $R$ and $R'$ are identical since we have not changed any arrow used by $R$ anyway. The sequences $L$ and $L'$ agree up to time $T$, while $L'_n\ge x$ for all $n\ge T$, since $L'$ can never go left from $x$ after time $T$.  It follows that $n_{L'}(z)=n_{L,T}(z)<\infty$ for every $z<x$.   

Let $y_1:=\max\{z<x:n_{R'}(z)\ge n_{L',T}(z)\}$.  By Lemma \ref{lem:exists}, $-\infty<y_1<x$.
By Lemma \ref{lem:x-1,x} (applied to $L'$, $R'$) either $R'$ hits $y_1+1$ at least $n_{L'}(y_1+1) \ge n_{L,T}(y_1+1)$ times, or $R'$ is forever to the right of $y_1+1$ after fewer than $n_{L'}(y_1+1)$ visits. In either case, $y_1+1<x$ (as $n_{R'}(x)<k$ and $R'$ lies eventually to the left of $x$). So there exists some $y_2\in (y_1,x)$ such that $n_{R'}(y_2)\ge n_{L'}(y_2)=n_{L',T}(y_2)$. This contradicts the definition of $y_1$.
\Qed

\begin{COR}
\label{cor:neighbour}
If $n_{R,t}(x-1)>n_{L,t}(x-1)$ then $n_{R,t}(x)\ge n_{L,t}(x)$.
\end{COR}
\proof Suppose instead that $n_{R,t}(x)< n_{L,t}(x)$.  Let $k=n_{R,t}(x)+1$, so that $T=T_L(x,k)\le t$ and $S=T_R(x,k)>t$.  Then 
\[n_{R,S}(x-1)\ge n_{R,t}(x-1)>n_{L,t}(x-1)\ge n_{L,T}(x-1).\]
This violates Lemma \ref{lem:main1} (if $n_R(x)\ge k$) or Lemma \ref{lem:main2} (if $n_R(x)<k$).
\Qed

\begin{COR}
\label{cor:right_hit_time}
Fix $x>0$, and let $T=T_{L}(x,1)=\inf\{t:L_t=x\}$ and $S=T_{R}(x,1)$.  Then $S\le T$.
\end{COR}
\proof If $T=\infty$ then the result is trivial.  So assume $T<\infty$. Lemma \ref{lem:main2} with $k=1$ implies that $S<\infty$ as well ($R$ cannot be to the right of $x>0$ without ever passing through $x$).  For each $i<x$, the number of times that $L$ hits $i$ before $T$ is $n_{L,T}(i)$, so $T=\sum_{i=-\infty}^{x-1}n_{L,T}(i)$.  Moreover, $n_{L,T}(i)$ is the number of times that $L$ hits $i$ before hitting $i+1$ for the $n_{L,T}(i+1)$-th time  (by definition of $T$, the last visit to $i<x$ up to time $T$ occurs before the last visit to $i+1$ up to time $T$).  By Lemma \ref{lem:main1} with $k=1$ we get that $n_{R,S}(x-1)\le n_{L,T}(x-1)$. Set $k_0=1$.

Now apply Lemma \ref{lem:main1} with $x-1$ instead of $x$ and with $k_1=n_{R,S}(x-1)$ to get
%From what we have just shown, $n_L(x-1)\ge n_{L,T}(x-1)\ge k$. By definition of $k$, it is not true that $R$ lies to the right of $x-1$ after fewer than $k$ visits to $x-1$. Therefore
$$
n_{R,T_R(x-1,k_1)}(x-2)\le n_{L,T_L(x-1,k_1)}(x-2).
$$
But $n_{R,T_R(x-1,k_1)}(x-2)=n_{R,S}(x-2)$ since $R$ cannot visit $x-2$ at times in $(T_r(x-1,k_1), S]$ (in other words, the last visit to $x-2$ occurs before the last visit to $x-1$). Furthermore, $n_{L,T_L(x-1,k_1)}(x-2)\le n_{L,T}(x-2)$ since $n_{L,T}(x-1)\ge k_1\Rightarrow T_L(x-1,k_1)\le T$. We have just shown that 
$$
n_{R,S}(x-2)=n_{R,T_R(x-1,k_1)}(x-2)\le n_{L,T_L(x-1,k_1)}(x-2)\le n_{L,T}(x-2).$$
Iterating this argument while $k_j=n_{R,S}(x-j)>0$ by applying Lemma \ref{lem:main1} at $x-j$ with $k=k_j$ (there is nothing to do once $n_{R,S}(x-j)=0$ for some $j$), we obtain by induction that $n_{R,S}(i)\le n_{L,T}(i)$ for every $i<x$. 
Thus $S=\sum_{i=-\infty}^{x-1}n_{R,S}(i)\le \sum_{i=-\infty}^{x-1}n_{L,T}(i)= T$ as required.\Qed

It follows immediately from Corollary \ref{cor:right_hit_time} that 
\eq
\overline{R}_n:=\max_{k\le n}R_k\ge \max_{k\le n}L_k=:\overline{L}_n.
\label{eq:maxbound}
\en 
Of course by mirror symmetry we also have $\underline{R}_n:=\min_{k\le n}R_k\ge \min_{k\le n}L_k=\underline{L}_n$. The following result extends this idea to the number of visits of the two paths to $\overline{R}_n$ by time $n$.

\begin{LEM}
\label{lem:maxvisits}
For each $t\ge0$, $n_{R,t}(\overline{R}_t)\ge n_{L,t}(\overline{R}_t)$ and $n_{L,t}(\underline{L}_t)\ge n_{R,t}(\underline{L}_t)$.
\end{LEM}
\proof Let $\mc{L}\CL \mc{R}$ and suppose the first claim fails.  Let $T=\inf\{t\ge 0: n_{R,t}(\overline{R}_t)< n_{L,t}(\overline{R}_t)\}<\infty$. 
Let $\mc{N}_t=n_{L,t}(\overline{R}_t)-n_{R,t}(\overline{R}_t)$.  Then $\mc{N}_{t+1}-\mc{N}_{t}\le 1$ if $\overline{R}_{t+1}=\overline{R}_t$, and by 
(\ref{eq:maxbound}),
$\mc{N}_{t+1}=0$ or $-1$ if $\overline{R}_{t+1}>\overline{R}_t$. Therefore by definition of $T$ we must have $R_T<\overline{R}_T$, $L_T=\overline{R}_T$, and $n_{L,T}(\overline{R}_T)=1+n_{R,T}(\overline{R}_T)$.
%\[
%\[n_{L,T}(\overline{R}_T)-1=n_{L,T-1}(\overline{R}_T)=n_{R,T-1}(\overline{R}_T)=n_{R,T}(\overline{R}_T).\]
%By (\ref{eq:maxbound}), $\overline{L}_{T-1}\le \overline{R}_{T-1}$, 
%so if $R_T=\overline{R}_{T-1}+1$ then $n_{R,T}(\overline{R}_T)=1$ and $n_{L,T}(\overline{R}_T)\le 1$. This would contradict the definition of $T$. Therefore in fact $R_T\le\overline{R}_{T-1}$ and so $\overline{R}_{T}=\overline{R}_{T-1}$. In particular, $n_{R,T-1}(\overline{R}_T)\ge n_{L,T-1}(\overline{R}_T)$. 
%To have $n_{R,T}(\overline{R}_T)<n_{L,T}(\overline{R}_T)$ we must have equality above, as well as $L_T=\overline{R}_T>R_T$. In this case, $n_{L,T}(\overline{R}_T)=1+n_{R,T}(\overline{R}_T)$. 
Moreover this happens regardless of the arrows of $\mc{L}$ or $\mc{R}$ at $\overline{R}_T$ above level $n_{R,T}(\overline{R}_T)$.  Define new arrow systems $\mc{L}',\mc{R}'$ by setting all arrows at $\overline{R}_T$ at level $1+n_{R,T}(\overline{R}_T)$ and above to be $\rightarrow$.  By construction $\mc{L}'\CL \mc{R}'$, and  $(L_n,R_n)=(L_n',R_n')$ for $n\le T$.  However $\overline{L}_{T+1}'=\overline{R}_T+1>\overline{R}_T=\overline{R}_{T+1}'$ which violates the fact that $\overline{R}'_n\ge \overline{L}'_n$ for all $n\ge 0$.   

The second result follows by mirror symmetry.\qed
\medskip

For each $z\in \Z$, $t \in \Z_+$, let $\overline{z}_t=\max(n_{L,t}(z),n_{R,t}(z))$.
\begin{LEM}
\label{lem:plus_minus_counts}
If there exist $t,y$ such that $R_t\le y<L_t$ and $n_{R,t}(y)>n_{L,t}(y)$ then $n_{R,t}(x)\ge n_{L,t}(x)$ for every $x\in [y,L_t]$.
\end{LEM}
\proof Suppose that $t$ and $y$ satisfy the above hypotheses, but the conclusion fails for some $x\in [y,L_t]$. In other words, $y<x\le L_t$ and $n_{R,t}(x)<n_{L,t}(x)$.  Define new arrow systems $\mc{L}'$ and $\mc{R}'$ by setting: 
\begin{itemize}
\item all arrows at $y$ at level $n_{R,t}(y)+I_{\{R_t\ne y\}}$ and above to be $\leftarrow $; 
\item all arrows at $x$ at level $n_{L,t}(x)+I_{\{L_t\ne x\}}$ and above to be $\rightarrow$; and 
\item for each $z>x$ set all arrows above level $\overline{z}_t$ to be $\rightarrow$.  
\end{itemize}
The resulting arrow systems satisfy $\mc{L}'\CL\mc{R}'$ with $(L_n,R_n)=(L_n',R_n')$ for $n\le t$.  By construction $L_n'\ra \infty $ as $n\ra \infty$, since $L_n'$ never again goes below $x$, and can make at most finitely many more $\leftarrow$ moves. But also  $R'_n\le y$ for all $n\ge t$, which contradicts the fact that $\overline{R}'_n\ge \overline{L}'_n$ for all $n\ge 0$.\qed
\medskip

We say that a sequence $\{L_n\}_{n\ge 0}$ on $\Z$ is {\em transient to the right} if for every $x\in \Z$ there exists $n_x\ge 0$ such that $L_n>x$ for all $n\ge n_x$ (i.e.~if $\liminf_{n\ra \infty}L_n=+\infty$). 
\begin{COR}
\label{cor:transience}
If $\liminf_{n\ra \infty}L_n=+\infty$ then $n_R(x)\le n_L(x)$ for every $x$ and $\liminf_{n\ra \infty}R_n=+\infty$.
\end{COR}
\proof  Suppose that $L$ is transient to the right.  Then $n_L(y)<\infty$ for each $y$.  Suppose that for some $x$, $n_R(x)>n_L(x)$.  Let $T=T_{R}(x,n_L(x)+1)$.  Define new systems $\mc{L}'\CL \mc{R}'$ by setting every arrow at $x$ above level $n_L(x)$ to be $\leftarrow$.  Then $L'=L$, so $L'\ra \infty$, but $R_t'\le x$ for every $t\ge T$.  This violates (\ref{eq:maxbound}) for $L'$, $R'$.  Therefore $n_R(x)\le n_L(x)$ for every $x$, which establishes the first claim.

For the second claim, suppose that $R$ is not transient to the right.  Then $R$ is either transient to the left or it visits some site $x$ infinitely often.  In either case there is some site $x$ such that $n_R(x)>n_L(x)$ which cannot happen by the first claim.
%the former case so that  $\overline{R}_{\infty}<\infty$ so (\ref{eq:maxbound}) is violated.  In the latter case, let $x$ be any site that $R$ visits  more times than $L$.  Define two new environments $\mc{L}'$ and $\mc{R}'$ by changing all arrows at site $x$ from level $n_L(x)+1$ and above to be $\leftarrow$.  Then we still have $\mc{L}'\CL \mc{R}'$ and $L'=L$, while $R'=R$ up to the time $T_{R}(x,n_L(x)+1)$, and $R'_n\le x$ for all $n\ge T_{R}(x,n_L(x)+1)$.  This violates (\ref{eq:maxbound}) for $L'$, $R'$.
\Qed

\begin{COR}
\label{cor:greater_i_o}
$R\ge L$ infinitely often.
\end{COR}
\proof If $R$ is not bounded above, this follows by considering the times at which $R$ extends its maximum. It follows similarly if $L$ is not bounded below, using times at which $L$ extends its minimum. The only remaining possibility is that $R$ is bounded above and $L$ is bounded below, in which case by (\ref{eq:maxbound}) 
both paths visit only finitely many vertices.  In this case consider the sets of vertices that $R$ and $L$ visit infinitely often.  
Let $x_{\infty}=\sup\{z\in \Z:n_{R}(z)=\infty\}$ and 
$y_{\infty}=\sup\{z\in \Z:n_{L}(z)=\infty\}$.
If $x_\infty<y_\infty$ then Lemma \ref{lem:main2} is violated (apply it to $x=y_\infty$ for $k>n_{R}(y_{\infty})$). Therefore $x_\infty\ge y_\infty$, so $R_t\ge L_t$ at all sufficiently large $t$ for which $R_t=x_\infty$. 
\Qed
%\begin{COR}
%\label{cor:last visit}
%Suppose that both $L$ and $R$ are transient to the right.  Let $x>0$ and let $T_x^U$ and $T_x^L$ be the times of the last visits to $x$ by $R$ and $L$ respectively.  Then $T_x^U\le T_x^L$.
%\end{COR}
%\proof Suppose that $L$ visits $x$ exactly $k$ times.  Then $R$ visits $x$ at most $k$ times (proved as in Lemma \ref{lem:main}).  Let $i\le k$ be the number of visits of $R$ to $x$.  Suppose that $r^+_U$ of these are from the right and 
\subsubsection{Proof of Theorem \ref{thm:main}}
To prove (i) we show that if $L_n\ge x$ for all $n$ sufficiently large, then $R_n\ge x$ for all $n$ sufficiently large.  Suppose instead that $R_n<x$ infinitely often.  Then choose $N$ sufficiently large so that $L_n\ge x$ for all $n\ge N$, but $R_N<x$ and $n_{R,N}(R_N)>n_{L,N}(R_N)$.  Define two new arrow systems $\mc{L}',\mc{R}'$ by switching all arrows at $R_N$ from level $n_{R,N}(R_N)$ and above to be $\leftarrow$.  Then $\mc{L}'\CL\mc{R}'$ but Lemma \ref{lem:plus_minus_counts} is violated, as is Corollary \ref{cor:greater_i_o}.  This establishes (i).  Applying (i) to $-\mc{R}\CL -\mc{L}$ establishes (ii).
%Similarly if $L_n\ge x$ infinitely often but $R_n\ge x$ only finitely often then there exists some $z\ge x$ such that $n_L(z)>n_R(z)$.  Let $S=\inf\{t:R_n<x \,\forall n\ge t\}$ and let $T=\inf\{t\ge S:n_{L,t}(z)=n_R(z)+1\}$.  Define new environments $\mc{L}^*\CL\mc{R}^*$ by switching every arrow at $z$ above level $n_R(z)$ to $\rightarrow$ and for every $y>z$ switching every arrow at $y$ above level $n_{L,T}(y)\vee n_{R}(y)$ to $\rightarrow$.  The resulting walks $L^*$ and $R^*$ eventually violate (\ref{eq:maxbound}) since $L_n^*\ra \infty$ and $\overline{R}^*_{\infty}<\infty$.  
%that violate Lemma \ref{lem:plus_minus_counts}.  
%This establishes the first two claims.

If $R_n\ge x$ infinitely often then  $\limsup R_n/a_n\ge \limsup x/a_n=0$.  Thus the result is trivial unless there exists $0< M<\infty$ such that $\limsup L_n/a_n>M$.
%If equality holds then the required result is trivial since we've just proved that $R_n\ge x$ infinitely often.  Otherwise let $0< M<\infty$ be such that $\limsup L_n/a_n>M$. 
Then $L_n$ visits infinitely many sites $>0$.
Let $T_i$ be the times at which $L$ extends its maximum, i.e.~$T_0=0$ and for $i\ge 1$, $T_i=\inf\{n>T_{i-1}:L_{n}=1+\max_{k<n}L_k\}$.  
We first verify the (intuitively obvious) statement that $\frac{L_{T_i}}{a_{T_i}}>M$ infinitely often.  If $\frac{L_{T_i}}{a_{T_i}}>M$ only finitely often then for all $i$ sufficiently large, $\frac{L_{T_i}}{a_{T_i}}\le M$.  But for all $n\in [T_i,T_{i+1})$, $\frac{L_n}{a_n}\le\frac{ L_{T_i}}{a_n}\le \frac{L_{T_i}}{a_{T_i}}$.  So $\frac{L_n}{a_n}\le M$ for all but finitely many $n$, contradicting the fact that $\limsup L_n/a_n>M$.  Let $S_i$ be the times at which $R$ extends its max.  By definition, $L_{T_i}=i=R_{S_i}$ and from Corollary \ref{cor:right_hit_time}, $i\le S_i\le T_i$.  It follows immediately that for infinitely many $i$,
\eq
\frac{R_{S_i}}{a_{S_i}}\ge \frac{L_{T_i}}{a_{T_i}}>M,\nn
\en
whence $\limsup_{n \ra \infty}\frac{R_n}{a_n}\ge M$. This establishes part (iii)

To prove (iv), suppose that (iv) does not hold, and let $\tau$ be the first time at which this fails.  In other words
\[\tau=\inf\{t\ge 0: \text{ there exist } y, x<y \text{ such that } n_{R,t}(x)>n_{L,t}(x) \text{ and }n_{R,t}(y)< n_{L,t}(y)\}.\]
  Let $x_0$ be the largest such $x$, i.e.~$x_0=\sup\{x\in \Z: n_{R,\tau}(x)>n_{L,\tau}(x), \exists y>x \text{ such that }n_{R,\tau}(y)<n_{L,\tau}(y)\}$ and $y_0=\inf\{y>x_0:n_{R,\tau}(y)<n_{L,\tau}(y)\}$.  Then $x_0\le y_0-2$ or else Corollary \ref{cor:neighbour} is violated.  
By definition of $x_0$ and $y_0$ we have $n_{R,\tau}(y_0-1)\ge n_{L,\tau}(y_0-1)$.  Let $k=n_{L,\tau}(y_0)$.  Then $n_{L,\tau}(y_0-1)\ge n_{L,T_L(y_0,k)}(y_0-1)$ so $n_{R,\tau}(y_0-1)\ge n_{L,T_L(y_0,k)}(y_0-1)$.  On the other hand $n_{R,\tau}(y_0)<k$, so $\tau <T_R(y_0,k)$.  If $R_{\tau}<y_0-1$ then $n_{R,T_R(y_0,k)}(y_0-1)\ge n_{R,\tau}(y_0-1)+1>n_{L,T_L(y_0,k)}(y_0-1)$.  This contradicts one of the Lemmas \ref{lem:main1} or \ref{lem:main2} (depending on whether $n_R(y_0)\ge k$), so we must have instead that $R_\tau\ge y_0-1>x_0$.  Therefore $n_{R,\tau-1}(x_0)=n_{R,\tau}(x_0)>n_{L,\tau}(x_0)\ge n_{L,\tau-1}(x_0)$. Similarly if $L_\tau>x_0+1$ we get a contradiction to the symmetric versions of Lemmas \ref{lem:main1} or \ref{lem:main2}, so we must have $L_\tau\le x_0+1<y_0$, and therefore $n_{L,\tau-1}(y_0)=n_{L,\tau}(y_0)>n_{R,\tau-1}(y_0)$.  This contradicts the definition of $\tau$.

Finally, to prove (v), note that if $\mc{L}\CL \mc{R}$ then also $\mc{L}_+\CL \mc{R}_+$.  If $\mc{R}$ is $0$-right recurrent, then $R_{+,n}=0$ infinitely often so $L_{+,n}=0$ infinitely often by (i).
%Therefore  $\underline{L}_{+,t}=0$ for every $t$, so by Lemma \ref{lem:maxvisits}, $n_{L_+,t}(0)\ge n_{R_+,t}(0)$. If $\mc{R}$ is $0$-right recurrent, then $n_{R_+,t}(0)\to\infty$, so $n_{L_+,t}(0)\to\infty$ as well, establishing (i).
\Qed

%\subsubsection{Proof of Theorem \ref{thm:recurrent}}  
%Let $|\mc{I}|\CL\mc{E}$. That is, $\mc{I}_+\CL\mc{E}_+$ and $\mc{E}_-\CL\mc{I}_-$. Therefore $\underline{I}_{+,t}=0$ %for every $t$, so by Lemma \ref{lem:maxvisits}, $n_{I_+,t}(o)\ge n_{E_+,t}(o)$. If $\mc{E}$ is $o$-right recurrent, %then $n_{E_+,t}(o)\to\infty$, so $n_{I_+,t}(o)\to\infty$ as well, establishing (i).  Part (ii) follows by symmetry.
%%Similarly, $\overline{I}_{-,t}=0$ for every $t$, so $n_{I_-,t}(o)\ge n_{E_-,t}(o)\to\infty$ if $\mc{E}$ is $o$-left %recurrent, establishing the second part of the theorem.
%\Qed

%$n_{R,\tau}(z)=n_{L,\tau}(z)$ for all $x_0<z<y_0$.  

%Let $k=n_{L,\tau}(y_0)$.  Then 
%%$n_{L,\tau}(y_0)>n_{R,\tau}(y_0)$, so that 
%$T_R(y_0,k)>\tau\ge T_L(y_0,k)$ but 
%\[n_{R,T_R(y_0,k)}(y_0-1)\ge n_{R,\tau}(y_0-1)\ge n_{L,\tau}(y_0-1)\ge n_{L,T_L(y_0,k)}(y_0-1),\]
%which violates Lemma \ref{lem:main}.
%%If $R_\tau\ne x_0$ then $L_\tau=y_0$ by definition of $\tau$.  Suppose without loss of generality that $L_{\tau}-y_0$.  If $R_\tau\le x_0$ then Lemma \ref{lem:plus_minus_counts} is violated.  If $x_0<R_{\tau}<y_0$ then 
%\Qed}

\section{Counterexamples}
\label{sec:counterex}
\subsection{$L\TL R$ does not imply that $\liminf\frac{L_n}{n}\le \liminf\frac{R_n}{n}$}
\label{sec:counterex1}
In general, $L\TL R$ does not imply that $\liminf\frac{L_n}{n}\le \liminf\frac{R_n}{n}$, as we shall see in the following example.  

Let us first define the two systems as follows, starting with $\mc{L}$.  At $0$ the first three arrows are $\ra$.  At every $x>0$ the first two arrows are $\leftarrow$ and the next three arrows are $\rightarrow$.  It is easy to check that such a system results in a sequence $L$ that takes steps with the pattern $\rightarrow \leftarrow \rightarrow \leftarrow \rightarrow$ repeated indefinitely (without ever needing to look at arrows other than those specified above).  Thus $\lim_{n\ra \infty}\frac{L_n}{n}=\frac{3-2}{5}=\frac{1}{5}$.  

\begin{figure}
\includegraphics[scale=.5]{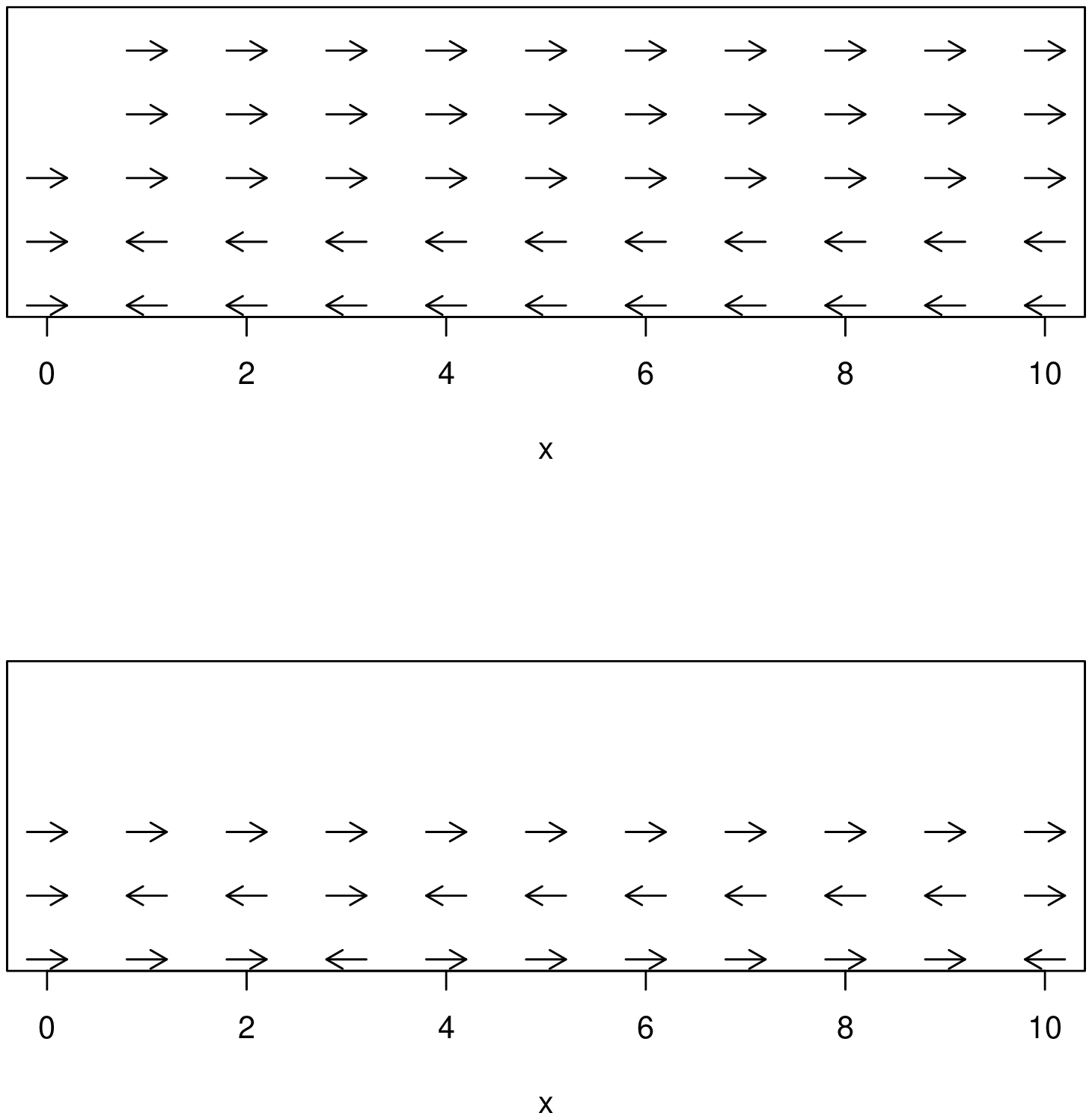}
\includegraphics[scale=.5]{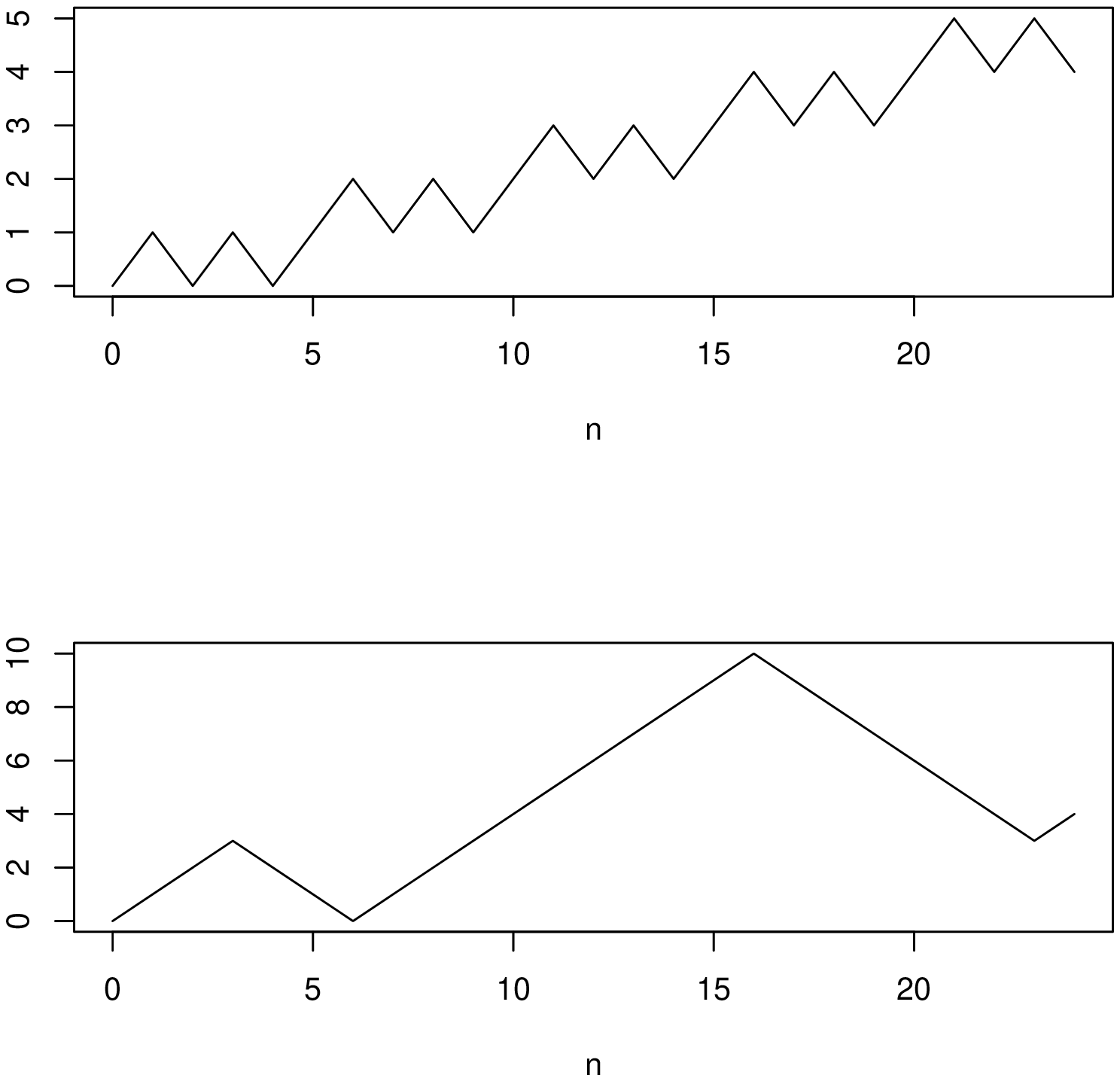}
\caption{On the left are parts of the systems $\mc{L}$ (top) and $\mc{R}(3)$ (bottom) and on the right are their corresponding sequences $R_n$ (solid) and $L_n$ (dotted), defined in Section \ref{sec:counterex} such that $\liminf n^{-1}L_n\ge \liminf n^{-1}R_n$.  Each site in $\N$ appears five times in the sequence $L$ and three times in the sequence $R$.}
\label{fig:counter}
\end{figure}

Let us now define a system $\mc{R}=\mc{R}(N)$, according to a parameter $N$ as follows.  At $0$ the first three arrows are $\ra$.  At each site $x_k=x_k(N)$ of the form 
\eq
x_k=\sum_{m=1}^kN^m-\sum_{m=1}^{k-1}\sum_{r=0}^m(-1)^{m-r}N^r, \qquad k\ge 1\label{eq:CEsites}
\en
the first arrow is $\leftarrow$ and the next two arrows are $\rightarrow$.   At all remaining sites $x>0$, the first three arrows are $\rightarrow,\leftarrow,\rightarrow$.   See Figure \ref{fig:counter} for parts of the systems $\mc{L}$ and $\mc{R}(3)$.  By definition of these systems the arrows to the left of 0 and above those shown are irrelevant, so we can set them to be the same (for example, all $\rightarrow$).

%par(mfrow=c(2,1))
%piece=c(0,1,0,1,0)
%plot(0:24,c(piece,piece+1,piece+2,piece+3,piece+4),xlab="n",ylab="",main="",type="l")
%Rn=c(0:3,2,1,0:10,9:3,4)
%plot(0:24,Rn,xlab="n",ylab="",main="",type="l")

%increments=c(1,-1,1,-1,1)
%allinc=rep(increments,4)
%Lwalk=c(0,cumsum(allinc))
%plot(0:length(allinc),Lwalk,type="l")

By construction $L\TL R$ for each $N\ge 1$, but we will show that $\liminf\frac{R_n}{n}\le \frac{1}{2N+1}<\frac{1}{5}$ for $N\ge 3$ (also $\limsup\frac{R_n}{n}\ge \frac{N}{N+2}$).

The first site of the form (\ref{eq:CEsites}) is $x_1=N$.  The walk $R$ first encounters a $\leftarrow$ at its first visit to this site and then sees a $\rightarrow$ at site 0 (second visit to 0). The walk $R$ then visits site $x_1$ for the second time, whence it sees a $\rightarrow$.  It continues moving right, visiting every site between $x_1$ and $x_2$ exactly once before reaching $x_2$ at this point it sees a $\leftarrow$, moves to $x_2-1$ (for the second visit to that site) and continues seeing $\leftarrow$ at every site in $(x_1,x_2)$ until reaching $x_1$ for the third time.  It then sees $\rightarrow$ at every site in $[x_1,x_2)$ (third visit to each of those sites), but also at every site in $[x_2,x_3)$ (second visit to $x_3$ and first visit to each site in $(x_3,x_4)$).  Continuing in this way, the walk turns left at every $x_i$ on the first visit, and continues left (second visit at interior sites) until reaching $x_{i-1}$ for the third time, and then continues to go right until reaching $x_{i+1}$ for the first time.

%The distance between two sites of the form (\ref{eq:CEsites}) is 
%\eqalign
%x_{k+1}-x_{k}
%%&=
%%\sum_{m=1}^{k+1}N^m-\sum_{m=1}^{k}\sum_{r=0}^m(-1)^{m-r}N^r-\left(\sum_{m=1}^{k}N^m-\sum_{m=1}^{k-1}\sum_{r=0}^m(-1)%^{m-r}N^r\right)\nn\\
%%&=N^{k+1}+\sum_{r=0}^{k}(-1)^{k+1-r}N^r=\sum_{r=0}^{k+1}(-1)^{k+1-r}N^r
%%\nn\\
%%&=(-1)^{k+1}\left[\frac{(-N)^{k+2}-1}{-N-1}\right]
%=\frac{N^{k+2}+(-1)^{k+1}}{N+1}.
%\enalign
%The lengths of the $k$th up and down periods respectively are therefore
%\eq
%x_{k}-x_{k-2}
%%=x_{k}-x_{k-1}+x_{k-1}-x_{k-2}=\sum_{r=0}^{k}(-1)^{k-r}N^r+\sum_{r=0}^{k-1}(-1)^{k-1-r}N^r
%=N^k, \quad \text{ and }x_k-x_{k-1}=\sum_{r=0}^{k}(-1)^{k-r}N^r.
%\en
%%while the length of the $k$th down period is $x_k-x_{k-1}=\sum_{r=0}^{k}(-1)^{k-r}N^r$.

At time $t_k=\sum_{m=1}^kN^m + \sum_{m=1}^{k-1}\sum_{r=0}^m(-1)^{m-r}N^r$ the walk is at position $x_{k}=\sum_{m=1}^{k}N^m-\sum_{m=1}^{k-1}\sum_{r=0}^m(-1)^{m-r}N^r$ for the first time.  Simple calculations then give
%At these times we have
%\eq
%\frac{R_{t_k}}{t_k}=\frac{\sum_{m=1}^{k}N^m-\sum_{m=1}^{k-1}\sum_{r=0}^m(-1)^{m-r}N^r}{\sum_{m=1}^kN^m + %\sum_{m=1}^{k-1}\sum_{r=0}^m(-1)^{m-r}N^r}.
%\en
%Relatively simple calculations then give
%where
%\eqalign
%\sum_{m=1}^{k-1}\sum_{r=0}^m(-1)^{m-r}N^r&=\sum_{m=1}^{k-1}(-1)^m\left[\frac{(-N)^{m+1}-1}{-(N+1)}\right]=\frac{N}{N+1}\sum_{m=1}^{k-1}N^m+\frac{1}%{N+1}\sum_{m=1}^{k-1}(-1)^m\nn\\
%&=\frac{N}{N+1}\left(\frac{N^k-1}{N-1}-1\right)+\frac{1}{N+1}\left[\frac{(-1)^k-1}{-2}-1\right]\\
%&=\frac{N}{N+1}\left(\frac{N^k-1}{N-1}-1\right)-\frac{1}{2(N+1)}\left[1-(-1)^k\right],
%\enalign
%and 
%\eq
%\sum_{m=1}^{k}N^m=\frac{N^{k+1}-1}{N-1}-1.
%\en
%It follows that 
\eqalign
\lim_{k\ra \infty}\frac{R_{t_k}}{t_k}&=
%\lim_{k\ra \infty}\frac{\frac{N^{k+1}-1}{N-1}-\frac{N^{k+1}}{(N+1)(N-1)}}{\frac{N^{k+1}-1}{N-1}+\frac{N^{k+1}}{(N+1)(N-1)}}=\lim_{k\ra \infty}\frac{(N+1)(N^{k+1}-1)-N^{k+1}}{(N+1)(N^{k+1}-1)+N^{k+1}}=
\frac{N}{N+2},\nn
\enalign
which gives rise to the limit supremum claimed.  

Similarly at times $s_k=\sum_{m=1}^kN^m + \sum_{m=1}^{k}\sum_{r=0}^m(-1)^{m-r}N^r$ the walk is at position $x_{k-1}=\sum_{m=1}^{k}N^m-\sum_{m=1}^{k}\sum_{r=0}^m(-1)^{m-r}N^r$ for the last time.  After some simple calculations we obtain 
\eqalign
\lim_{n\ra \infty}\frac{R_{s_k}}{s_k}&=
%\lim_{k\ra \infty}\frac{\frac{N^{k+1}-1}{N-1}-\frac{N^{k+2}}{(N+1)(N-1)}}{\frac{N^{k+1}-1}{N-1}+\frac{N^{k+2}}{(N+1)(N-1)}}=\lim_{k\ra \infty}\frac{(N+1)(N^{k+1}-1)-N^{k+2}}{(N+1)(N^{k+1}-1)+N^{k+2}}=
\frac{1}{2N+1},\nn
\enalign
which gives rise to the limit infimum claimed.  
\subsection{$L$ can be in the lead more than $R$}
\label{sec:counterex2}
Given two sequences $L$ and $R$ with $L\CL R$, let $A_{R,t}=\{n\le t: R_n>L_n\}$ and $A_{L,t}=\{n\le t: R_n<L_n\}$.  It is not unreasonable to expect that for every $t\in \mathbb{N}$, $|A_{R,t}|\ge |A_{L,t}|$ which essentially says that $R$ is ahead of $L$ more than $L$ is ahead of $R$.  It turns out that this does not hold even when $L\TL R$.

\begin{figure}
\includegraphics[scale=.4]{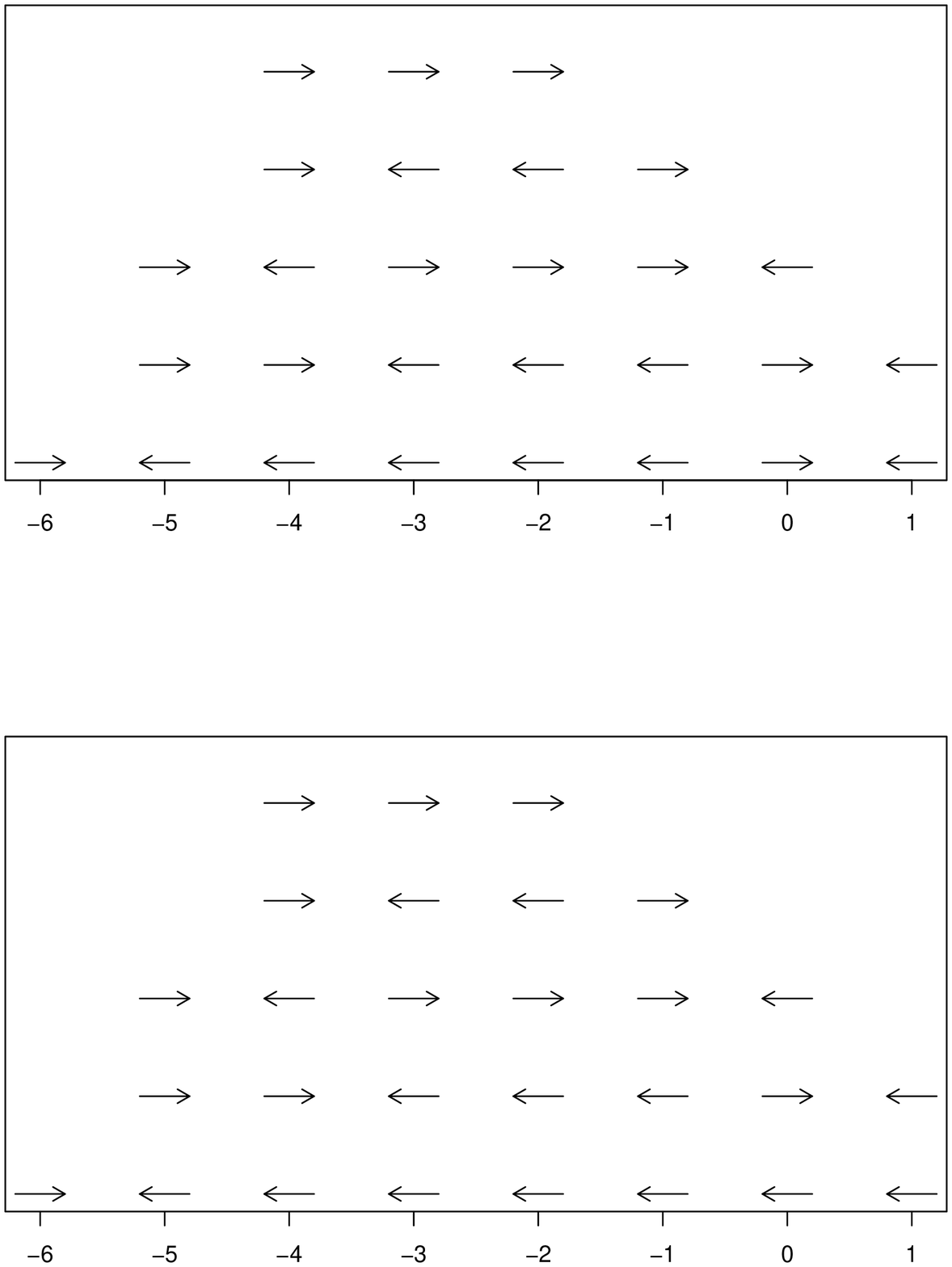}
\includegraphics[scale=.4]{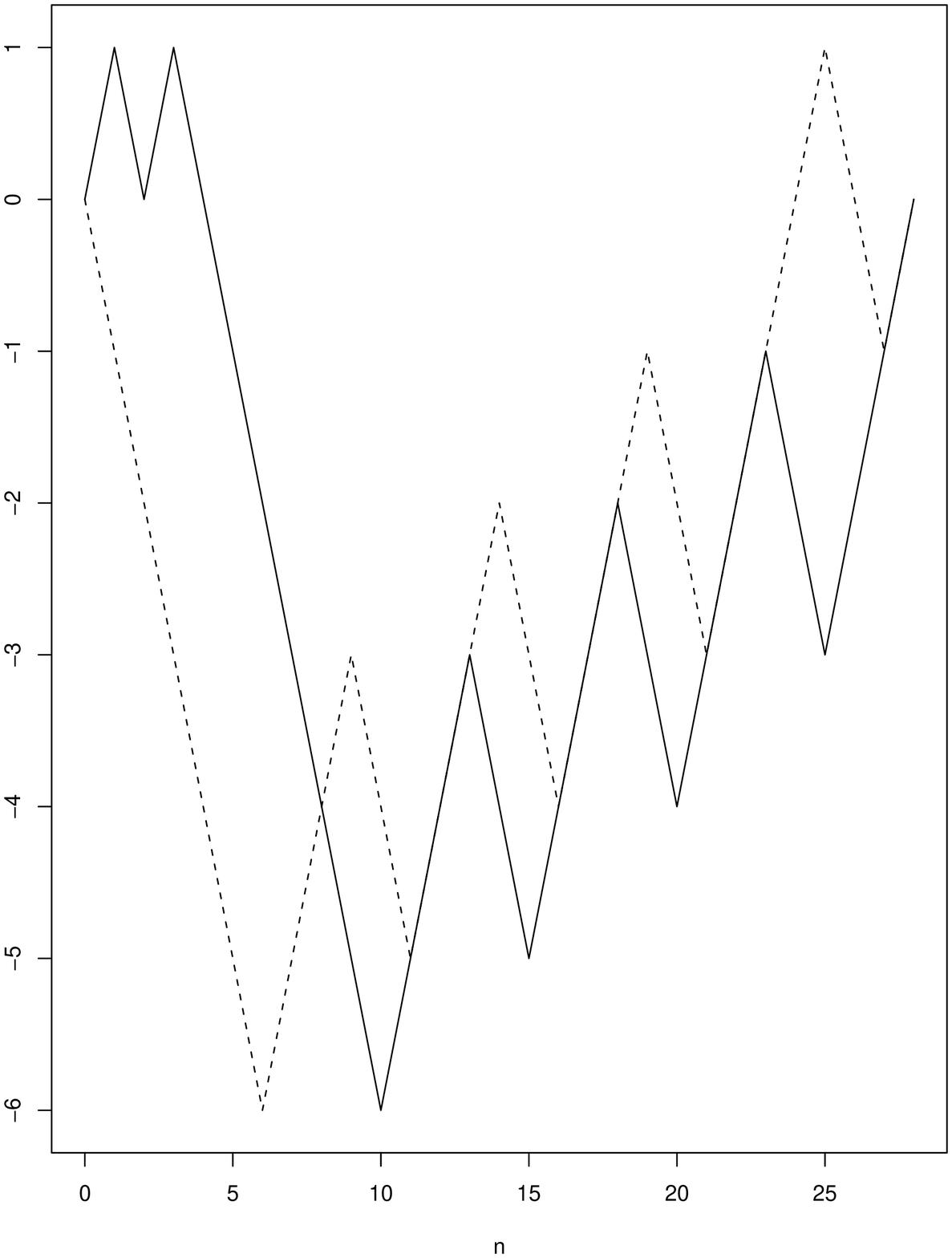}
\caption{On the left are parts of arrow systems $\mc{L}$ (top) and $\mc{R}$ (bottom) with $\mc{L}\TL \mc{R}$, and on the right are the corresponding paths $R_n$ (solid) and $L_n$ (dotted).  Here, $|A_{R,26}|=7<9= |A_{L,26}|$.}
\label{fig:leader}
\end{figure}

To see this, consider the partial arrow systems $\mc{R}$ and $\mc{L}$ on the left hand side of Figure \ref{fig:leader}.  These two systems differ only at the first arrow at $0$, whence $\mc{L}\TL\mc{R}$ (if we set all other arrows to be equal, for example).  The first 28 terms of the sequences $L$ and $R$ are plotted on the right of the figure.  At any place where the solid line is above the dotted line, $R>L$.  In particular $R_n>L_n$ only for $1\le n\le 7$.  Similarly $L>R$ when the dotted line lies above the solid line, which happens at times $9,10,14,15,19,20,24,25,26$.  Thus we have  $|A_{R,25}|=7<8= |A_{L,25}|$ and similarly $|A_{R,26}|=7<9= |A_{L,26}|$.

\begin{figure}
\begin{center}
\includegraphics[height=7cm,width=7cm]{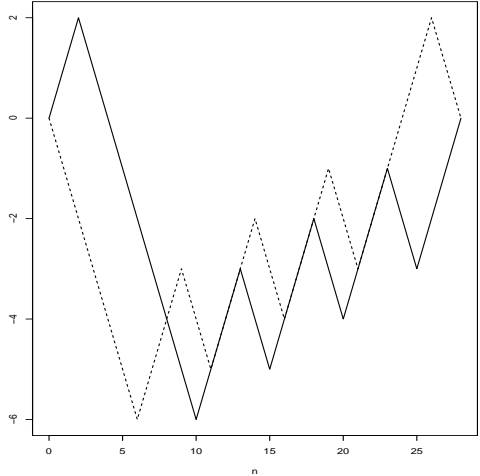}
\end{center}
\caption{Paths $R'_n$ (solid) and $L'_n$ (dotted) with $L_n'\CL R_n'$ and $|A_{R,28}|=7<10= |A_{L,28}|$.  The walks have visited each site the same number of times.}
\label{fig:leader2}
\end{figure}

We can modify these systems slightly to get another interesting example.  Define $\mc{R}'$ from $\mc{R}$ by switching the second arrow at 0 to $\leftarrow$, the first arrow at 1 to be $\rightarrow$ and setting the first arrow at 2 to be $\leftarrow$.  Define $\mc{L}'$ from $\mc{L}$ by switching the first arrow at 1 to be $\rightarrow$ and setting the first arrow at 2 to be $\leftarrow$.  The resulting partial systems satisfy $\mc{L}'\CL\mc{R}'$.  At time $t=28$, $|A_{R,28}|<|A_{L,28}|$, the number of visits to each site is identical, and $L_{28}=R_{28}=0$ (see Figure \ref{fig:leader2}).  This means we can define a system which repeats such a pattern indefinitely.   We can add any common steps that we wish in between repetitions of this pattern and hence we can have recurrent, transient, or even ballistic sequences satisfying $L\CL R$ but 
such that $t^{-1}(|A_{L,t}|-|A_{R,t}|)\ra v>0$ as $t \ra \infty$.

\section{Applications}
\label{sec:applications}
In this section we describe some of the applications of our main results in the theory of nearest neighbour self-interacting random walks, i.e.~sequences $(X_n)_{n\ge 0}$ of $\Z$-valued random variables (which may include projections of higher dimensional walks), such that $X_{n+1}-X_n\in \{-1,1\}$ a.s.~for every $n$.  For each application, what we actually do is show that there is a probability space on which the relevant random walks live and on which they are related via the property $\CL$ or $\TL$ almost surely.  It is then clear that on that probability space the conclusions of Theorem \ref{thm:main} hold almost surely for the walks satisfying those relations.

Our original motivation for the present paper was in studying random walks in (non-elliptic) random environments in dimensions $d\ge 2$ (see e.g.~\cite{HS_RWDRE}).  In \cite{HS_RWDRE} the authors apply Theorem \ref{thm:main} to  random walks in i.i.d.~random environments such that for some diagonal direction $u$, with sufficiently large probability at each site there is a drift in direction $u$, and that almost surely there is no drift in direction $-u$.  For such walks, the projection $R$ in direction $u$ can be coupled with a so-called 1-dimensional multi-excited random walk (see below) $L$ so that $L\TL R$, and transience and positive speed results can be obtained for this projection, when the strength of the drift is sufficiently large. 

Our results can also be applied to recurrent models.  For example, given $\beta>-1$, let $X$ be a once-reinforced random walk (ORRW) on $\Z$ with reinforcement parameter $\beta$, i.e.~$X_0=0$ and 
\eq
\label{eq:ORRW}
\mP(X_{n+1}-X_n=1|\mc{F}_n)=\frac{1+\beta I_{\{X_{n}+1\in \vec{X}_{n-1}\}}}{2+\beta [I_{\{X_{n}+1\in \vec{X}_{n-1}\}}+I_{\{X_{n}-1\in \vec{X}_{n-1}\}}]}.\nn
\en
We can similarly define ORRW on $\Z^+$ by forcing the walk to step right when at $0$.  Then it is possible to define a probability space on which there is a ORRW $X^+(\beta)$ for each $\beta>-1$ and such that $X^+(\beta)\TL X^+(\zeta)$ whenever $\beta\ge \zeta>-1$.  On this probability space the corresponding local times processes then satisfy the monotonicity property Theorem \ref{thm:main}(iv).

Most of our results, including that for random walks in random environments above, involve comparisons with so-called {\em multi-excited random walks in i.i.d.~cookie environments}.  A cookie environment is an element ${\bf \omega}=(\omega(x,n))_{x \in \Z,n \in \mathbb{N}}$ of $[0,1]^{\Z\times \mathbb{N}}$.  A (multi-)excited random walk in cookie environment ${\bf \omega}$, starting from the origin, is a sequence of random variables $X=\{X_n\}_{n\ge 0}$ defined on a probability space (and adapted to a filtration $\mathcal{F}_n$) such that $X_0=0$ a.s.~and
\[P_{\bf \omega}\big(X_{n+1}=X_n+1\big|\mc{F}_n\big)=\omega(x,\ell(n))=1-P_{\bf \omega}\big(X_{n+1}=X_n-1\big|\mc{F}_n\big),\]
where $\ell(n)=\ell_X(n)=\sum_{m=0}^n1_{\{X_m=X_n\}}$.  In other words, if you are currently at $x$ and this is the $k$th time that you have been at $x$ then your next step is to the right with probability $\omega(x,k)$, independent of all other information.  A random cookie environment ${\bf \omega}$ is said to be i.i.d.~if the random vectors $\omega(x,\cdot)$ are i.i.d. as $x$ varies over $\mathbb{Z}$. 

Let ${\bf U}=(U(x,n))_{x \in \Z,n \in \mathbb{N}}$ be a collection of independent standard uniform random variables defined on some probability space.  For each $x\in \Z$, $n\in \N$, and each cookie environment ${\bf \omega}$ let
\eq
\mc{E}_{{\bf \omega},{\bf U}}(x,n)=\begin{cases}
\rightarrow, & \text{ if }U(x,n)<\omega(x,n)\\
\leftarrow, & \text{ otherwise.}\label{eq:randomL}\nn
\end{cases}
\en
Then $\mc{E}_{{\bf \omega},{\bf U}}$ is an arrow system determined entirely by the pairs $(\omega(x,n),U(x,n))_{x\in \Z,n \in \N}$, and the corresponding walk $E=E_{{\bf \omega},{\bf U}}$ is an excited random walk in cookie environment ${\bf \omega}$.  Given two cookie environments ${\bf \omega}$ and ${\bf \omega'}$ we write ${\bf \omega}\TL {\bf \omega'}$ if $\omega(x,n)\le \omega'(x,n)$ for every $x\in \Z$ and $n\in \N$.  If ${\bf \omega}\TL {\bf \omega'}$, then on the above probability space $\mc{E}_{{\bf \omega},{\bf U}}\TL \mc{E}_{{\bf \omega'},{\bf U}}$ so Theorem \ref{thm:main} applies to the corresponding excited random walks.

%\eq
%\mc{R}_{{\bf \omega'},{\bf U}}(x,n)=\begin{cases}
%\rightarrow, & \text{ if }U(x,n)<\omega'(x,n)\\
%\leftarrow, & \text{ otherwise,}
%\end{cases}
%\en
For excited random walks in i.i.d.~cookie environments in 1 dimension, it is known up to a high level of generality that right transience and the existence of a positive speed $v>0$ do not depend on the order of the cookies (see e.g.~\cite{KZ08}).  One might expect that the value of $v$ should depend on this order.  
%It is therefore natural to wonder whether or not there is an $\CL$ analogue of Theorem \ref{thm:mono1}.  The answer is yes.  
The main result of this section is Theorem \ref{thm:cookie_cl} below, which essentially states that one cannot decrease the ($\limsup$)-speed of a cookie random walk by swapping stronger cookies in a pile with weaker cookies that appear earlier in the same pile (and doing this at each site).  In order to state the result precisely we require some further notation.

For each $x\in \Z$, let $\mc{A}_x$ denote a partition of $\N$ into finite (non-empty) subsets.  For any such partition we can order the elements of the partition as $\mc{A}_x=(A^1_x,A^2_x,\dots)$ (e.g.~according to the ordering of the smallest element in each $A^i_x$).  Let $\mc{A}=(\mc{A}_x)_{x\in \Z}$ denote a particular collection of such partitions (indexed by $\Z$), and $\mc{P}$ denote the set of all such collections.  Let $\mc{P}_n$ denote the set of such collections where every $A_x^s$ is a set containing at most $n$ elements.  

Fix $\mc{A}\in \mc{P}$.  Let $x\in \Z$, $s\in \N$, $\omega$ be a cookie-environment, and $j,k\in A^s_x$ with $j\le k$.  We say that $(j,k)$ is an {\em $(x,s,\omega)$-favourable swap} if $\omega(x,j)\le \omega(x,k)$.  
Let $\omega(x,A^s_x)=(\omega(x,r))_{r\in A^s_x}$, and let $b=(j,k)$ be an $(x,s,\omega)$-favourable swap.  Define $\omega_b(s,A^s_x)$ by,
\[\omega_b(s,r)=\begin{cases}
\omega(s,k), &\text{ if }r=j\\
\omega(s,j), &\text{ if }r=k\\
\omega(s,r), &\text{ if }r\in A^s_x\setminus\{j,k\}.\end{cases}\] 
Then we say that $\omega_b(s,A^s_x)$ is the $A^s_x$-environment produced by the swap $b=(j,k)$, and write $\omega(x,A^s_x)\cb \omega_b(x,A^s_x)$.  Given two cookie environments ${\bf \omega}$ and ${\bf \omega'}$, we say that ${\bf \omega'}$ is {\em an $\mc{A}$-permutation of ${\bf \omega}$} if for each $s$ and $x$, $\omega'(x,A^s_x)$ is a permutation of $\omega(x,A^s_x)$.  If ${\bf \omega'}$ is an $\mc{A}$-permutation of ${\bf \omega}$ and 
%is {\em $\mc{A}$-favoured over $\omega$}
 if also on every $A^s_x$, ${\bf \omega'}$ can be generated from ${\bf \omega}$ from a finite sequence of favourable swaps then we write ${\bf \omega}\ACL{\bf \omega'}$.  More precisely ${\bf \omega}\ACL{\bf \omega'}$ if for every $x\in \Z$, $s\in \N$, $j\le k$, there exists a finite sequence of pairs of $A_x^s$ indices $b_1,\dots,b_K$ (for some $K\ge 0$), and $A^s_x$-environments $(\omega_i(x,A^s_x))_{i=0}^K$ with $\omega_0(x,A^s_x)=\omega(x,A^s_x)$ and $\omega_K(x,A^s_x)=\omega'(x,A^s_x)$ such that $\omega_i(x,A^s_x)\overset{b_{i+1}}{\rightarrow} \omega_{i+1}(x,A^s_x)$ are favourable swaps for each $i=0,\dots,K-1$. 

Given $\mc{A}\in \mc{P}$ and an environment ${\bf \omega}$, let $\underline{\bf \omega}_{\mc{A}}$ denote the environment obtained by permuting $\omega$ on each $A_x^s$ so that $\underline{\bf \omega}_{\mc{A}}(x,j)\le \underline{\bf \omega}_{\mc{A}}(x,k)$ for all $j,k \in A_x^s$ such that $j<k$.  Note that $\underline{\omega}_{\mathcal{A}}(x,A_x^s)$ can be obtained from $\omega(x,A_x^s)$ by a sequence consisting of at most $|A_x^s|-1$ swaps that are not favourable: first perform the swap that moves the largest $\omega(x,k)$ for $k\in A_x^s$ to the highest location in $A_x^s$, then proceed iteratively, always moving the next largest value to the next highest location.  Reversing this procedure generates $\omega(x,A_x^s)$ from $\underline{\omega}_{\mathcal{A}}(x,A_x^s)$ by a sequence of (at most $|A_x^s|-1$) favourable swaps, so that $\underline{\omega}_{\mathcal{A}} \preccurlyeq^{\mathcal{A}}\omega$. 

%To see this, let $A_x^s=(k_1,\dots,k_n)$ with $k_i< k_{i+1}$ for each $i$.  Let $k_{(1)},\dots, k_{(n)}$ be an $\omega$-ranking of the elements of $A_x^s$, so that $\omega(x,k_{(i)})\le \omega(x,k_{(i+1)})$ for each $i$. Starting from $\omega(x,A_x^s)$, perform the (unfavourable) swap $(k_{(n)},k_n)$ that moves the largest $\omega(x,k)$ for $k\in A_x^s$ to the highest location in $A_x^s$.  Proceed iteratively, moving the next largest $\omega(x,k)$ to the second highest position.  

%let $k_1 > k_2>\dots$ be the elements of $A^s_x$, move the appropriate element of
%$\underline{\omega}_{\mathcal{A}}$ to align with $\omega(x,k_1)$ via a sequence of favourable %swaps, then do more favourable swaps to align another element with  $\omega(x,k_2)$, etc.

For fixed ${\bf \omega}$ and for any finite subset $A_x\subset \N$, let $\{V_i\}_{i \in A_x}$ be a collection of i.i.d.~standard uniform random variables and define $N_{A_x}=\sum_{i\in A_x}I_{V_i\le \omega(x,i)}$ (which can be thought of as the number of right arrows generated by $\omega(x,A_x)$).  Note that the law of $N_{A_x}({\bf \omega})$ is invariant under permutations of the indices in the set $A_x$, so that $q_{{\bf \omega},A_x}(y)=\mP(N_{A_x}({\bf \omega})=y)$ is invariant under such permutations.

\begin{THM}
\label{thm:cookie_cl}
Let $\mc{A}\in \mc{P}_3$ and let ${\bf \omega}$ be a cookie environment.  Then there exists a probability space on which: for each $\mc{A}$-permutation ${\bf \omega'}$ of $\underline{\bf \omega}_{\mc{A}}$ there is an excited random walk $E_{\bf \omega'}$ in environment ${\bf \omega'}$, defined such that $E_{\bf \omega'}\CL E_{\bf \omega''}$ almost surely whenever ${\bf \omega'}\ACL{\bf \omega''}$.
\end{THM}
\proof Let ${\bf U}=\{U_{x,s}\}_{x\in \Z,s\in \N}$ be i.i.d.~standard uniform random variables, and ${\bf Y}=\{Y_{x,s}\}_{x\in \Z,s\in \N}$ be independent random variables (independent of ${\bf U}$) where $Y_{x,s}$ has the law of $N_{A_x^s}({\bf \omega})$ for each $x,s$.

Let $x\in \Z$ and $s\in \N$ and consider the set $A_x^s$, which contains $n=|A_x^s|\le 3$ elements.  Without loss of generality let us assume that $A_x^s=\{1,\dots,n\}$.  Let $y=Y_{x,s}$ and note that (since $n\le 3$) the set $S_{n,y}$ of $n$-stacks  (an $n$-stack is any element of $\{\leftarrow,\rightarrow\}^n$) containing exactly $y$ right arrows is a completely ordered set (under $\CL$) of cardinality $n_y={n \choose y}$.  Let $(a_1^{(y)},\dots,  a_{n_y}^{(y)})$ be the reverse ordering of the set (so that $a_1^{(y)}$ is the element 
%$(\rightarrow, \rightarrow, \dots, \rightarrow, \leftarrow, \leftarrow,\dots,  \leftarrow)$ 
consisting of $y$ right arrows underneath $n-y$ left arrows), and let $a_i^{(y)}(j)$ be the $j$th arrow of $a_i^{(y)}$.  

Now for any $\mc{A}$-permutation ${\bf \omega'}$ of  $\underline{\bf \omega}_{\mc{A}}$, define a probability measure $P_{{\bf \omega'}}$ on $S_{n,y}$ by setting 
\[P_{{\bf \omega'}}(a_i^{(y)}) =\big(q_{{\bf \omega},A_x^s}(y)\big)^{-1}\prod_{j=1}^n\left[\omega'(x,j)I_{a_i^{(y)}(j)=\rightarrow}+(1-\omega'(x,j))I_{a_i^{(y)}(j)=\leftarrow}\right], \quad i=1,\dots, n_y.\]
This is the conditional probability of selecting (for the arrows corresponding to $A_x^s$) a particular configuration $a_i^{(y)}$ consisting of $y$ right arrows and $n-y$ left arrows, given that the configuration contains exactly $y$ right arrows and $n-y$ left arrows. 
Define $\mc{E}_{{\bf \omega'}}(x,A_x^s)=(\mc{E}_{{\bf \omega'}}(x,j))_{j \in A_x^s}$ by 
\[\mc{E}_{{\bf \omega'}}(x,A_x^s)=a_m^{(y)}, \quad  \text{ if }\sum_{i=1}^{m-1}P_{{\bf \omega'}}(a_i^{(y)})<U_{x,s}\le \sum_{i=1}^{m}P_{{\bf \omega'}}(a_i^{(y)}).\]

Let ${\bf \omega'}$ and ${\bf \omega''}$ be $\mc{A}$-permutations of $\underline{\bf \omega}_{\mc{A}}$ with ${\bf \omega'}\ACL{\bf \omega''}$.  Recall that $q_{{\bf \omega''},A_x^s}(y)=q_{{\bf \omega'},A_x^s}(y)$ by invariance under permutations.  Also note that for every $m\le n_y$,
\[\sum_{i=1}^{m}P_{{\bf \omega''}}(a_i^{(y)})\ge \sum_{i=1}^{m}P_{{\bf \omega'}}(a_i^{(y)}),\]
so that under this coupling, $\mc{E}_{{\bf \omega'}}(x,A_x^s)=a_m^{(y)}\Rightarrow \mc{E}_{{\bf \omega''}}(x,A_x^s)=a_k^{(y)}$ for some $k\le m$.  This means that $\mc{E}_{{\bf \omega'}}(x,A_x^s)\CL \mc{E}_{{\bf \omega''}}(x,A_x^s)$ when we consider $\CL$ on $A_x^s$ only.

Let us now summarize what we have achieved.  For fixed $\mc{A}$ and ${\bf \omega}$, we have coupled arrow systems (and hence the corresponding walks) defined from all $\mc{A}$-permutations of $\underline{\bf \omega}_{\mc{A}}$ (including ${\bf \omega}$ itself) so that $\mc{E}_{{\bf \omega'}}(x,A_x^s)\CL \mc{E}_{{\bf \omega''}}(x,A_x^s)$ for each $x\in \Z$, $s\in \N$ when ${\bf \omega'}\ACL {\bf \omega''}$, where the coupling took place independently (according to the variables ${\bf U}$ and ${\bf Y}$) for each $x,s$.  It follows that for any such ${\bf \omega'}$, ${\bf \omega''}$, under this coupling, $\mc{E}_{{\bf \omega'}}\CL \mc{E}_{{\bf \omega''}}$.  The result follows since for each $\mc{A}$ permutation ${\bf \omega'}$, the corresponding walk $E_{{\bf \omega'}}$ has the law of an excited random walk in cookie environment ${\bf \omega'}$.
\Qed

Note that in the statement (and proof) of Theorem \ref{thm:cookie_cl} the probability space depends on $\mc{A}$ and ${\bf \omega}$ and is constructed in such a way that each $A^s_x$ has the same number of right arrows under ${\bf \omega}$ as under ${\bf \omega}'$ (and likewise left arrows).  If $\mc{A}\in \mc{P}_2$, which corresponds to considering only disjoint transpositions/swaps, then the above proof can be simplified slightly, and the probability space defined independently of ${\bf \omega}$).  The coupling is then defined on $\mc{A}_x^s=(j,k)$ for each ${\bf \omega}$ by
\eq
(\mc{E}(x,j),\mc{E}(x,k))=\begin{cases}
(\rightarrow,\rightarrow), & \text{ if }U_{x,k,j}<\omega(x,j)\omega(x,k)\\
(\rightarrow,\leftarrow), & \text{ if }\omega(x,j)\omega(x,k)\le U_{x,k,j}<\omega(x,j)\\
(\leftarrow,\rightarrow), & \text{ if }\omega(x,j)\le U_{x,k,j}<\omega(x,j)+\omega(x,k)(1-\omega(x,j))\\
(\leftarrow,\leftarrow), & \text{ otherwise.}\label{eq:randomL2}
\end{cases}
\en
This works because the set of 2-stacks is totally ordered according to $\CL$ as 
\[\begin{matrix}\rightarrow\\ \rightarrow \end{matrix}\quad \succcurlyeq \quad 
\begin{matrix}\leftarrow\\ \rightarrow \end{matrix}\quad \succcurlyeq \quad 
\begin{matrix}\rightarrow\\ \leftarrow \end{matrix}\quad \succcurlyeq \quad 
\begin{matrix}\leftarrow\\ \leftarrow \end{matrix}\]
so there is no need to define the random variables $Y_{x,s}$ whose laws depend on ${\bf \omega}$.  If on the other hand we relax the condition that $\mc{A}\in \mc{P}_3$ to $\mc{A}\in \mc{P}_4$ the proof breaks down because e.g.~the $4$-stacks $\RLLR$ and $\LRRL$ are not ordered by $\CL$.  However, by considering finite sequences of favourable swaps, we can obtain the following theorem.
\begin{THM}
\label{thm:cookie_cl2}
Let $\mc{A}\in \mc{P}$ and let ${\bf \omega}\ACL {\bf \omega'}$ be two cookie environments.  Then there exists a probability space on which there are excited random walks $E_{\bf \omega}$ and $E_{\bf \omega'}$ in environments ${\bf \omega}$ and ${\bf \omega'}$ respectively, defined such that $E_{\bf \omega}\CL E_{\bf \omega'}$ almost surely.
\end{THM}
\proof 
Fix $x\in \Z$, $s\in \N$.  Then ${\bf \omega'}(x,A_x^s)$ can be obtained from ${\bf \omega}(x,A_x^s)$ by a finite sequence of favourable swaps $\omega_i(x,A^s_x)\overset{b_{i+1}}{\rightarrow} \omega_{i+1}(x,A^s_x)$, $i=0,\dots, K_x^s-1$, with $\omega_0(x,A^s_x)=\omega(x,A^s_x)$ and $\omega_{K_x^s}(x,A^s_x)=\omega'(x,A^s_x)$.  Using the coupling in Theorem \ref{thm:cookie_cl} for a single favourable swap on $A_x^s$, for each $i$ we can define a probability space with finite chunks of random arrow systems $(\mc{E}_i(x,A_x^s),\mc{E}_i'(x,A_x^s))$ with marginal laws defined by $\omega_i(x,A^s_x)$ and $\omega_{i+1}(x,A^s_x)$ respectively, and such that $\mc{E}_i(x,A_x^s)\CL \mc{E}_i'(x,A_x^s)$.  

Let $(X_1,Y_1)$ and $(Y_2,Z_2)$ be random quantities (not necessarily defined on the same probability space) such that $Y_1$ and $Y_2$ have the same distribution.  Then we can construct $X_3$, $Y_3$, and $Z_3$ on a common probability space by letting $Y_3\sim Y_1\sim Y_2$, and letting $X_3$ and $Z_3$ be conditionally independent given $Y_3$, with conditional laws the same as $X_1$ given $Y_1$ and $Z_2$ given $Y_2$ respectively.  Iterating this construction, and applying the resulting coupling to the random objects $\mc{E}_i(x,A_x^s)$, we can construct a probability space on which there are finite chunks of random arrow systems $\mc{E}_i(x,A_x^s)$ with marginal laws defined by $\omega_i(x,A^s_x)$, $i=0,\dots, K_x^s$, such that $\mc{E}_i(x,A_x^s) \CL \mc{E}_{i+1}(x,A_x^s)$ for each $i$.  Taking the product probability space over $x\in \Z$ and $s\in \N$, and letting $\mc{E}=(\mc{E}_0(x,A_x^s))_{x\in \Z, s\in \N}$ and $\mc{E}'=(\mc{E}_{K_x^s}(x,A_x^s))_{x\in \Z, s\in \N}$, we have that $\mc{E}\CL \mc{E}'$.  Defining $E_{\omega}$ and $E_{\omega'}$ to be the corresponding walks gives the result.
\Qed

%However, one can of course recover some stochastic inequalities about random walks in environments ${\bf \omega}\ACL{\bf \omega'}$ for $\mc{A}\in \mc{P}_n$ (for any $n$) by using Theorem \ref{thm:cookie_cl} iteratively to define a finite sequence of probability spaces and performing just one favourable swap per $A^s_x$ at each iteration.  We believe that for any $n$, $\mc{A}\in \mc{P}_n$ and ${\bf \omega}\ACL{\bf \omega'}$ it is possible to couple corresponding arrow systems $\mc{E}_{\bf \omega}\CL \mc{E}_{\bf \omega'}$ on a single probability space

Since Theorems \ref{thm:cookie_cl} and \ref{thm:cookie_cl2} are defined rather abstractly, we now give an explicit example.  
Suppose that ${\bf \omega}$ is an environment defined by $\omega(x,2k-1)=p_1$ and $\omega(x,2k)=p_2$ for every $x\in \Z$, $k\in \N$, with $p_2>p_1$.  Suppose also that we wish to understand the effect (on the asymptotic properties of the corresponding excited random walk) of switching the order of the first two cookies at every even site, or instead, of switching the values of $p_1$ and $p_2$ at even sites.  In the first case the environment of interest is ${\bf \omega'}$ where $\omega'(x,1)=\omega(x,2)$ and  $\omega'(x,2)=\omega(x,1)$ for each $x\in 2\Z$ and otherwise $\omega'(x,k)=\omega(x,k)$, while in the second case we have ${\bf \omega''}$ defined by $\omega''(x,2k-1)=\omega(x,2k)$ and $\omega''(x,2k)=\omega(x,2k-1)$ for all $x\in 2\Z$, $k\in \N$ and otherwise $\omega'(x,k)=\omega(x,k)$.  In this example the permutations of interest are composed of {\em disjoint} swaps/transpositions, and hence we can choose partitions consisting of sets containing at most 2 elements.  For example, letting $A_x^s=\{2s-1,2s\}$ for each $x\in \Z$, $s\in \N$ defines one particular choice (among many) of $\mc{A}$ for which ${\bf \omega'}$ and ${\bf \omega''}$ are $\mc{A}$-permutations of ${\bf \omega}$, and such that ${\bf \omega}\ACL {\bf \omega'}\ACL {\bf \omega''}$.  Theorems \ref{thm:cookie_cl} and \ref{thm:main} then imply that e.g.~if $p_1\ge \hlf$ (so that the walks are not transient to the left) then the limsup speeds of the corresponding random walks satisfy $\overline{v}_{\omega}\le \overline{v}_{\omega'}\le \overline{v}_{\omega''}$.

The ORRW is an example of a walk whose drift can depend on more than just the number of visits to the current site.  For example, on $\Z_+$ the drift encountered by the ORRW at a site $x$ at time $n$ (so $X_n=x$) depends on whether the local time of the walk at $x+1$ is positive.  Some of the known results for excited random walks in i.i.d.~or ergodic environments can be extended to more general self-interacting random walks (where the drifts may depend on the history in an unusual way) with a bounded number of positive drifts per site.
\begin{THM}
\label{thm:general}
Let $X_n$ be a nearest-neighbour self-interacting random walk and $\mc{F}_n=\sigma(X_k,k\le n)$.  Suppose that there exist $M\in \mathbb{N}$ and $(\eta_k)_{k\le M}\in [0,1)^M$ such that
\begin{itemize}
\item $\mP(X_{n+1}=X_n+1|\mc{F}_n)I_{\ell(n)=k}\le \eta_k$ for all $k\le M$ and all $n\in \Z_+$ almost surely, and 
\item $\mP(X_{n+1}=X_n+1|\mc{F}_n)I_{\ell(n)=k}\le \hlf$ for all $k> M$ and all $n\in \Z_+$, almost surely.
\end{itemize}
If $\alpha=\sum_{k=1}^M(2\eta_k-1)\le 1$ then $X$ is not transient to the right, almost surely.  If $\alpha\le 2$ then $\limsup n^{-1}X_n\le 0$, almost surely.  If $\alpha<-1$ then $X$ is transient to the left, almost surely.  If $\alpha <-2$ then $\liminf n^{-1}X_n<0$ almost surely.
\end{THM}

\proof  Define $\eta_k=\hlf $ for $k>M$.  For each $x\in \Z$, let $\omega(x,k)=\eta_k$ for $k \in \N$.  Let ${\bf U}=(U(x,m))_{x\in \Z,m\in \N}$ be i.i.d.~standard uniform random variables.   and
define $\mc{R}$ by 
\[\mc{R}(x,k)=\begin{cases}
\rightarrow &\text{ if }U(x,k)\le \eta_k\\
\leftarrow &\text{ otherwise.}
              \end{cases}\]
The corresponding walk $R_n$ has the law of an excited random walk in the (non-random) environment $\omega$.  By \cite{KZ08}, the conclusions of the theorem hold for the walk $R$, e.g.~if $\alpha=\sum_{k=1}^M(2\eta_k-1)\le 1$ then $R$ is not transient to the right, almost surely.  

For a nearest neighbour sequence $x_0,\dots,x_n$ define 
$$
P_{n,k}(x_0,\dots,x_n)=\mP(X_{n+1}=X_n+1|X_0=x_0, \dots, X_n=x_n)I_{\ell_{x}(n)=k}.
$$ 
Define a nearest neighbour self-interacting random walk $L$ by setting $L_0=0$ and given that $\ell_{L}(n)=k$,
\[L_{n+1}=\begin{cases}
L_n+1 , & \text{ if } U(L_n,k)\le P_{n,k}(L_0,\dots,L_n)\\
L_n-1, &\text{ otherwise.}          \end{cases}\]
Then $L$ has the law of $X$.  Since $P_{n,k}\le \eta_k$ almost-surely, we have that $L\TL R$ almost surely.  The result now follows by Cor. \ref{cor:transience}.
The astute reader may have noticed that we have not defined the arrow system $\mc{L}$.  We can do so, according to the walk $L$ as follows.  Given that $\ell_{L}(n)=k$, define 
\[\mc{L}(L_n,k)=\begin{cases}
\rightarrow, &\text{ if }U(L_n,k)\le P_{n,k}(L_0,\dots,L_n)\\
\leftarrow, &\text{ otherwise}.
\end{cases}\]
In other words, this inductively defines $\mc{L}$ as the arrow system determined by the steps of the walk $L$.  Since $L$ does not define an entire arrow system at any site $x$ visited only finitely often by $L$ we can define $\mc{L}(x,k)=\leftarrow$ for each $k>n_{L}(x)$.  

To be more precise, for each $n$ we can define $\mc{L}^{(n)}$ according to the arrow system determined by $L_0,\dots,L_n$ and adding $\leftarrow$ everywhere else.  For each such $n$ we have $\mc{L}^{(n)}\TL \mc{R}$, so that Theorem \ref{thm:main} (iv) holds for each $n$, and so does (\ref{eq:maxbound}).  The former result implies the claims about transience when $\alpha\le 1$ and $\alpha<-1$, while (\ref{eq:maxbound}) and its minimum equivalent imply the remaining results (see e.g.~the proof of Theorem \ref{thm:main} (iii)).
\Qed

\section*{Acknowledgements}
The authors would like to thank two anonymous referees for their helpful suggestions and the Fields Institute for hosting them while part of this work was carried out.  This research was supported in part by the Marsden Fund (Holmes), and by NSERC (Salisbury).

%\bibliography{C:/Users/mhol003/Desktop/texfiles/maaak}
\bibliographystyle{plain}

\end{document}